\documentclass[english,12pt,leqno]{article}
\usepackage{tikz}
\usetikzlibrary{decorations.markings}
\usetikzlibrary{calc}

\usepackage[T2A]{fontenc}

\usepackage[latin1, russian, utf8]{inputenc}

\usepackage{amsxtra}
\usepackage{amsmath}
\usepackage{amssymb}
\usepackage{amsfonts}
\usepackage[all]{xy}
\usepackage{mathrsfs}
\usepackage{amsthm}
\usepackage{float}
 \usepackage{caption}

\newtheorem{thm}[equation]{Theorem}

\newtheorem{cor}[equation]{Corollary}
\newtheorem{lem}[equation]{Lemma}
\newtheorem{prop}[equation]{Proposition}

\newtheoremstyle{example}{\topsep}{\topsep}%
     {}
     {}
     {\bfseries}
     {.}
     {2pt}
     {\thmname{#1}\thmnumber{ #2}\thmnote{ #3}}

   \theoremstyle{example}
   
   \newtheorem{Defi}[equation]{Definition}
   
   \newtheorem{rem}[equation]{Remark}
   \newtheorem{rems}[equation]{Remarks}

   \newtheorem{ex}[equation]{Example}

 \newtheorem{exer}[equation]{Exercise}

\newtheoremstyle{example}{\topsep}{\topsep}%
     {}
     {}
     {\bfseries}
     {.}
     {2pt}
     {\thmname{#1}\thmnumber{ #2}\thmnote{ #3}}

   \numberwithin{equation}{section}

\def\CC{\mathbb{C}}

\def\FF{\mathbb{F}}

\def\PP{\mathbb{P}}
\def\RR{\mathbb{R}}
\def\ZZ{\mathbb{Z}}

\def\HH{\mathbb{H}}

\def\gen{\mathfrak{g}}

\def\hen{\mathfrak{h}}
\def\len{\mathfrak{l}}

\def\Sen{\mathfrak{S}}
\def\sen{\mathfrak{s}}

\def\Ac{\mathcal{A}}
\def\Bc{\mathcal{B}}
\def\Cc{\mathcal{C}}
\def\Kc{\mathcal{K}}
\def\Dc{\mathcal{D}}
\def\Ec{\mathcal{E}}
\def\Fc{\mathcal{F}}
\def\Gc{\mathcal{G}}

\def\Mc{\mathcal{M}}

\def\Oc{\mathcal{O}}
\def\Pc{\mathcal{P}}

\def\Rc{\mathcal{R}}
\def\Sc{\mathcal{S}}

\def\Vc{\mathcal{V}}

\def\Sen{\mathfrak{S}}

\def\fb{\mathbf{f}}

\def\fb{\mathbf{f}}

\def\<{\langle}
\def\>{\rangle}

\def\be{\begin{equation}}
\def\ee{\end{equation}}

\def\bef{\begin{figure}[H]\centering}
\def\enf{\end{figure}}
\def\Br{\on{Br}}

\def\btp{\begin{tikzpicture}}
\def\etp{\end{tikzpicture}}

\def\coh{{\on{coh}}}

\def\Cone{\on{Cone}}

\def\Ch{\on{Ch}}

\def\ds{{\Delta^{\on{sim}}}}

\def\Hom{\on{Hom}}
\def\hra{\hookrightarrow}

\def\Id{{\on{Id}}}

\def\k {\mathbf k}
\def\Ker{\on{Ker}}

\def\LG{\on{LG}}

\def\lla{\longleftarrow}
\def\lra{\longrightarrow}

\def\max{{\on{max}}}
\def\Mod{\on{Mod}}

\def\lra{\longrightarrow}

\def\ol{\overline}
\def\on{\operatorname}

\def\Perv{\on{Perv}}
\def\pt{{\on{pt}}}

\def\Ss{\mathbb S}
\def\Sch{\on{Sch}}

\def\ul{\underline}

\def\Vect{\on{Vect}}

\title{ Perverse Schobers}

\author{ Mikhail Kapranov, Vadim Schechtman}

\begin{document}


 \maketitle

 
 The notion of  a perverse sheaf, introduced in \cite{BBD}, has come to play a central role in algebraic
 geometry and representation theory. In particular, appropriate categories of perverse sheaves provide
 ``categorifications''  of various representation spaces, these spaces being recovered as the Grothedieck groups
 of the categories. 

  The goal of this paper is to suggest the possibility of categorifying the very concept of a perverse sheaf.
   In other words,  we propose to develop a theory of perverse sheaves not of vector spaces but of
 {\em triangulated  categories}. 
 
 \vskip .2cm
  
 Given a complex manifold  $X$, an analytic Whitney stratification $\Sc= (X_\alpha)_{\alpha\in A}$ of $X$ and
 a ground field $\k$, one has the category $\Perv(X,\Sc)$ of  perverse sheaves of $\k$-vector
 spaces on $X$ smooth with respect to $S$.  Traditionally, there have been two ways of looking at $\Perv(X,\Sc)$:
 
 \begin{itemize}
 \item[(1)] {\bf General definition:} as an abelian subcategory in the triangulated category $D^b_{\on{constr}}(X,\Sc)$
 of constructible complexes of sheaves of $\k$-vector spaces on $X$, smooth with respect to $\Sc$.
 
 \item[(2)] {\bf Quiver description} (for some particular $(X,\Sc)$){\bf :}  as a category of diagrams of 
 some given type
  formed by  vector spaces $(V_i)_{ i\in I}$ and maps between them subject
 to certain relations. These diagrams have the following features:
 
 \begin{itemize}
 
 \item[(2a)] Arrows  come in pairs 
 $
 \xymatrix{
 V_i \ar@<.4ex>[r]&V_j \ar@<.4ex>[l]
 }
 $
 having the same ends but opposite directions. 
 This reflects the (Verdier) self-duality of $\Perv(X,\Sc)$. 
 
 \item[(2b)] In most cases, the relations contain 2 or 3 summands, with coefficients $\pm 1$. 
 
 \end{itemize}

 \end{itemize}
 
 So far, there is no obvious
 direct way to categorify the approach (1) since it is not clear what are complexes of triangulated categories. 
 
 \vskip .2cm
 
 On the other hand, we observe that the features (2a) and (2b) {\em  are of the kind that immediately
 suggest a categorical generalization. } We can replace vector spaces $V_i$ by triangulated categories
 $\Vc_i$ and arrows by exact functors. The pairs of opposite arrows in (2a) can be interpreted
 as adjoint pairs of functors, 2-term relations as isomorphisms of functors and 3-term ones as exact
 triangles in  appropriate functor categories (see Appendix for the
 precise framework in which these make sense). 
 By forming Grothendieck groups\footnote{or applying any other functor from 
 the category of triangulated categories to $\k$-vector spaces, 
 for example, the higher Quillen $K$-theory, or Hochschild cohomology}
 $V_i = K_0(\Vc_i)\otimes\k$ of such a diagram of categories,
one would then obtain a quiver in the original sense, i.e., a perverse sheaf. 
The idea that the cone can be seen as a categorical analog of the difference, lies, of course,
at the very foundations of algebraic K-theory, especially in the Waldhausen approach. 

\vskip .2cm

This strongly suggests that  there should be meaningful objects which can be understood
as ``perverse sheaves of triangulated categories" and which give usual perverse sheaves
by passing to the Grothendieck groups.  We propose to call such hypothetical objects {\em perverse Schobers}
(or, sometimes, for brevity, simply {\em Schobers}), 
using the
German analog\footnote{
  A literal Russian analog would be the word
 стог. We learned the term ``Schober" from W. Soergel. 
}
 of the English word ``stack"  which would be the correct (but overused)   term for speaking of
``sheaves of categories". 

\vskip .2cm

In this paper we work out several basic examples of quiver descriptions of perverse sheaves
and define, in an {\em ad hoc} way,  what should be the perverse Schobers in these situations.
In the simplest case, we propose, in \S \ref{sec:disk}, to identify   perverse Schobers on a disk 
with one allowed singular point, 
with {\em spherical functors} of  \cite{A1}\cite{AL2}.

\vskip .2cm

 For a disk with several  allowed singular points we
propose, in \S \ref {sec:disk-several}, a definition in terms of certain diagrams of spherical functors,
and explain the invariance properties of such a definition. Among other things, we reformulate the
classical  Picard-Lefschetz
formula  as a general statement about perverse sheaves on a disk, and then lift it  to a distinguished
triangle associated to a perverse Schober and a certain configuration of paths. Such ``Picard-Lefschetz triangles"
should therefore be considered as fundamental features of perverse Schobers. It is natural to expect
analogous features in the case $\dim (X)>1$ and even consider them as
``codimension 1 data" of a perverse Schober. 

\vskip .2cm

A series of examples of spherical functors is provided by representation theory. More precisely,  for a
reductive group $G$ we have spherical
  functors acting in the derived category of sheaves on
$G/B$ and  satisfying  the relations  of the corresponding braid group $\Br(\gen)$. 
In \S \ref{sec:G/B} we review these examples and suggest a conjectural interpretation in terms of perverse Schobers on 
$\hen/W$.

\vskip .2cm

Most of the known quiver descriptions of $\Perv(X,\Sc)$
can be obtained using a choice of ``cuts" which are certain totally real subvarieties  $K\subset X$ (of real dimension equal to $\dim_\CC X$). 
In \S \ref{sec:fuk} we summarize the features of such cuts and note that Lagrangian varieties, used in constructing
Fukaya categories \cite{FOOO}\cite{FSbook},  provide a reasonable class of candidates for cuts. One can therefore expect Fukaya-categorical
constructions to have a bearing on the problem of classification of perverse sheaves. 

\vskip .2cm

Additionally, we discuss the idea of defining
``Fukaya categories with coefficients". This idea was proposed by M. Kontsevich in order to study the usual Fukaya category
of a manifold by fibering it over a manifold of smaller dimension. 
We suggest 
  that  perverse Schobers
should be considered as the right ``coefficient data" for such a definition, just like sheaves are natural coefficient
data for defining cohomology.

\vskip .2cm

We would like to thank A. Bondal, V. Ezhov, M. Finkelberg,  D. Nadler,
  P. Schapira,  W. Soergel,   Y. Soibelman and B. To\"en
for useful discussions and correspondence.
 V.S. is grateful to the Kavli IPMU for
hospitality and support during the visits when this paper was finished. 
The work of M. K.  was supported by World Premier International Research Center Initiative (WPI Initiative), MEXT, Japan.

\vfill\eject

   \section{Perverse Schobers on a disk: spherical functors }\label{sec:disk}
   
  \noindent {\bf A. Disk with one marked point.}  
  Let $\Delta$  be the unit disk
   in $\CC$ and   $\Perv(\Delta, 0)$ be the category of perverse sheaves on $\Delta$ with 
   the only possible singularity at $0$. The 
  most iconic example of a quiver description of perverse sheaves is  the  following  classical
  statement \cite{beil-gluing} \cite{GGM}.

   \begin{thm}\label{thm:ggm}
     $\Perv(\Delta,0)$ is equivalent to the category $\Pc_1$ of quadruples
   $(\Phi, \Psi, u, v)$ where $\Phi, \Psi\in \Vect_\k$ and
   \be\label{eq:phi-psi}
   \xymatrix{
 \Phi \ar@<.4ex>[r]^v&\Psi \ar@<.4ex>[l]^u
 }
   \ee
   are linear maps such that 
   \be\label{eq:t-psi}
   T_\Psi: = \Id_\Psi - vu \text{ is an isomorphism. }
   \ee
   \end{thm}
   
   \begin{exer}\label{exer:phi-psi}
 Show that $T_\Psi$ is an isomorphism iff 
$T_\Phi: \Id_\Phi - uv$ is an isomorphism. 
 \end{exer}  
   
   \bef
\centering
\btp[scale=.4, baseline=(current  bounding  box.center)]
\node (0) at (0,0){}; 
\fill (0) circle (0.15);
   
 \draw (0,0) circle (5cm);

   \node (b) at (5,0){}; 
\fill (b) circle (0.15);

\draw[line width = .3mm] (0,0) .. controls (3,.7) .. (5,0); 

\node at (-.7, .5) {$0$}; 
\node at (5.7, .5) {$b$}; 
\node at (2.7,1.3) {$K$}; 
\node at (4.5,4.5) {$\Delta$};

\etp

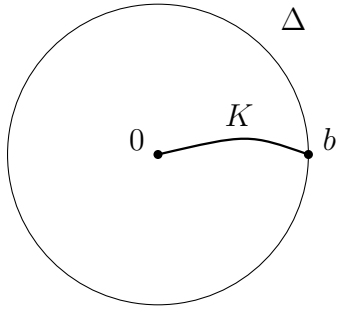
\caption{ Defining $\Phi$ and $\Psi$ topologically.  }
\label{fig:disk-1}
\enf
\noindent One way of constructing an explicit equivalence is as follows \cite{GGM}. Choose a base point $b$ on the boundary of $\Delta$ and
connect it with $0$ by a simple arc $K$, see Fig. \ref{fig:disk-1}. Then to $\Fc\in\Perv(\Delta, 0)$ we associate the spaces
\[
\begin{gathered}
\Phi(\Fc)\, \,:= \, \,\ul \HH^1_K(\Fc)_0  \,\,\simeq\,\,  \HH^1_K(\Delta, \Fc) \quad \text{(vanishing cycles)},\\
\Psi(\Fc)\,\, :=\,\, \Fc_b\,\, \simeq \,\, \ul \HH^1_K(\Fc)_b  \quad \text{(nearby cycles)}. 
\end{gathered}
\]
(We recall that $\Fc|_{\Delta-\{0\}}$ is a local system in degree $0$). 
The map $v=v_\Fc$ is the generalization map \cite{gelfand-macpherson}\cite{curry} 
for the constructible sheaf $\ul\HH^1_K(\Fc)$ on $K$, and $u$ is   the composition
\be\label{eq:map-u}
 \xymatrix{
 \Fc_b \ar[rr]^{\hskip -1cm
{ \text{counterclockwise} \atop \text{ continuation}  } }  && H^0(\Delta-K, \Fc)  \ar[r]^{\hskip .4cm \delta} & \HH^1_K(\Delta, \Fc).
 }
\ee

 \begin{rem}\label{rem:patt-1}
Further, we have the following elementary statements which we recall here as indicative of a certain pattern.
   \begin{enumerate}
\item[(1)]  $\ul\HH^i_K(\Fc) = 0$  { for }  $i\neq 1$.

\item[(2)] The sheaf $\Rc(\Fc) = \ul\HH^1_K(\Fc)$ on $K$ is constant on $K-\{0\}$ so it has only two essentially different stalks $\Rc(\Fc)_0 = \Phi(\Fc)$
and $\Rc(\Fc)_b=\Psi(\Fc)$. 

\item[(3)] Each of the two stalks, considered as a functor $\Perv(\Delta,0)\to\Vect_\k$, is an exact functor which takes
Verdier duality to vector space duality.

\item[(4)] The map   $u=u_\Fc$  is the dual
$u_\Fc = (v_{\Fc^*})^*$. 

\end{enumerate}

 \end{rem}

   \vskip .3cm
   
   \noindent {\bf B. Spherical functors.}
      As a natural categorical analog of the data \eqref{eq:phi-psi}-\eqref{eq:t-psi}  we would like to suggest  the following
  remarkable concept introduced by R. Anno \cite{A1}. 
  
  \begin{Defi}
  Let 
  \[
  S:\Dc_0 \lra \Dc_1
  \]
  be an exact functor between triangulated categories (see Appendix for conventions). Assume that 
 $S$ admits 
a left adjoint $L$ and a right adjoint $R$, so that we have the unit and counit natural transformations
\[
\begin{gathered}
SR\Rightarrow \Id_{\Dc_1}, \quad LS\Rightarrow\Id_{\Dc_0},\\
\Id_{\Dc_0}\Rightarrow RS, \quad \Id_{\Dc_1}\Rightarrow SL,
\end{gathered}
\] 
whose cones will be denoted by
\[
\begin{gathered}
T_1= \Cone\{SR\Rightarrow \Id_{\Dc_1}\}, \quad T'_1 = \Cone\{ \Id_{\Dc_1}\Rightarrow SL\} [-1] \,\,\,\text{ (the twist functors)},
\\
T_0 = \Cone \{ \Id_{\Dc_0}\Rightarrow RS \}[-1], \quad 
T'_0 = \Cone \{  LS\Rightarrow\Id_{\Dc_0}\}  \,\,\,\text{ (the cotwist functors).}
\end{gathered}
\]
We call $S$ a {\em spherical functor}, if:
\begin{enumerate}
\item[(SF1)] $T_1$ is an equivalence.

\item[(SF2)] The composition
$
R \to RSL\to T_0 L[1]
$
is an isomorphism. 
In other words, composition with $T_0$ identifies $R$ and $L$. 
\end{enumerate}
In this case $T'_1$ is quasi-inverse to $T_1$ and $T'_0$ is quasi-inverse to $T_0$. 
 \end{Defi}
 
 More precisely, each of the pairs of adjoint functors 
 \be\label{eq:schober-c-0}
  \xymatrix{
 \Dc_0 \ar@<.4ex>[r]^S&\Dc_1 \ar@<.4ex>[l]^R
 }, \quad 
  \xymatrix{
 \Dc_1\ar@<-.4ex>[r]_L&\Dc_0 \ar@<-.4ex>[l]_S
 }
 \ee
 can be regarded as an analog of \eqref{eq:phi-psi}.  Further,  (SF1) is an analog of \eqref{eq:t-psi}, 
 the adjunction unit allowing us to take the ``categorical difference", i.e., the cone. 
  Further still, the categorical analog of Exercise \ref{exer:phi-psi} can be found in the following result of
  Anno and Logvinenko \cite{AL2}.
  
  \begin{thm}\label{thm:AL}
  In addition to (SF1)  and (SF2) consider the following two conditions:
  \begin{enumerate}
  \item[(SF3)] $T_0$ is an equivalence.
  
  \item[(SF4)] The composition $L T_1[-1] \to LSR \to R$ is an isomoprhism. 
  \end{enumerate}
  Then, any two of the conditions (SF1)-(SF4) imply the other two. 
  \end{thm}
  
  \vskip .2cm
  
  So we can consider  a diagram \eqref{eq:schober-c-0}, i.e., the data of
  a spherical functor,  as the data defining  a ``perverse Schober''  over 
$(\Delta, 0)$. By passing to $K_0$ (or to any $K_i$, or to any homological functor, Hochschild homology for example) we get a perverse sheaf 
over $(\Delta, 0)$. 

\vskip .3cm

\noindent {\bf B. Examples of spherical functors.}  We now give some examples, to be used later.

\begin{ex}\label{ex:sph-p1}
Let $\Ss^d$ be the $d$-dimensional sphere and
 $q: \Ss^d \to \pt$ be the projection. We then have the functor
\[
\Dc_0 :=  D^b(\Vect) = D^b (\pt) \buildrel S=q^{-1}\over\lra D^b (\Ss^d) =:  \Dc_1
\]
with right adjoint  $R=Rq_*$ and left adjoint $L\simeq Rq_*[-d]$. The second adjunction is the Poincar\'e duality. Formaly, it
comes from the adjoint pair $(Rq_!, q^!)$ by noticing that $q^!\simeq  q^{-1}[d]$ (since $q$ is smooth orientable of relative
 dimension d) and that  $Rq_!=Rq^*$ (since $q$ is proper), see \cite{kashiwara-schapira}  for background. 
 More intrinsically,  $L$  is  canonically 
 identified with the tensor product of $Rq_*[-d]$ and $H^d(\Ss^d, \k)$,
  the 1-dimensional  $\k$-vector space spanned by global orientations of $\Ss^d$. 
 
 \end{ex}

\begin{prop}\label{prop:p1-sph}
 (a)  $S$ is a spherical functor. 
 
 (b) More generally, for any   $\Ss^d$-fibration $q: Z\to Y$ of  CW-complexes, the functor
 \[
 \Dc_0: = D^b(Y) \buildrel S= q^{-1}\over\lra D^b (Z) =: \Dc_1
 \]
 is a spherical functor. 
\end{prop} 

\noindent {\sl Proof:} We prove (a), since (b), being a relative version, is proved in the same way.

 The functor  $T_1=\Cone\{SR\Rightarrow\Id_{\Dc_1}\}$ is the endomorphism of  $D^b (S^d)$ defined as follows.
Let 
\[
U=(\Ss^d\times \Ss^d) - \Delta\,\, \buildrel j\over\hookrightarrow \,\,  \Ss^d \times \Ss^d
\]
be the embedding of the complement of the diagonal,  $\pi_1, \pi_2: \Ss^d \times \Ss^d \to \Ss^d $ be the projections
and $p_1, p_2: U\to \Ss^d$ be their restrictions to $U$.
Then
\[
T_1(\Fc) =  Rp_{2!}(p_1^*\Fc), 
\]
a formula remindful of the Fourier-Sato transform relating  sheaves on dual spheres \cite{SKK}. 
To see this, we write the functor in the RHS in terms of  a   ``kernel", as 
\[
 R\pi_{2*}((\pi_1^*\Fc)\otimes_\k \Kc), \quad \Kc = j_!\ul\k_U [1]\,\in\,  D^b (\Ss^d \times \Ss^d)
\] 
and note the exact sequence
\[
0\to \Kc[-1]\lra \ul\k_{\Ss^d \times \Ss^d } \lra \ul\k_\Delta\to 0,
\]
in which $\ul\k_{\Ss^d \times \Ss^d }$ is the kernel for $SR$ while $\ul\k_\Delta$ is the kernel for $\Id_{\Dc_1}$. 

 The twist $T'_1 = \Cone\{ {\Id_{\Dc_1}}\Rightarrow SL\}[-1]$ can be found explicitly as
\[
T'_1(\Fc) =  Rp_{2*}(p_2^*\Fc))[-1] = R\pi_{2*}((\pi_1^*\Fc)\otimes_\k \Kc'), \quad \Kc'= Rj_* \k_U [-1]. 
\]
Now, the condition (SF1), i.e., the fact that $T_1$ and $T'_1$ are quasi-inverse to each other, 
can be established directly by finding the ``composition'' of the kernels $\Kc$ and $\Kc'$
\[
\Kc * \Kc' \,\,=\,\, R\pi_{13*} (\pi_{12}^{-1}\Kc \otimes_\k \pi_{23}^{-1}\Kc'), \quad \pi_{ij}: \Ss^d\times\Ss^d\times\Ss^d \to \Ss^d\times\Ss^d,
\]
and showing that both $\Kc*\Kc'$ and $\Kc' * \Kc$ are isomorphic to $\ul \k_\Delta$ in degree 0. 
This amounts to the fact that for distinct points $x,y\in\Ss^d$ we have
\[
H^\bullet \left(\Ss^d  - \{x\}, j_! \bigl( \ul\k_{\Ss^d-\{x,y\}}\bigr) \right) \,=\, 0,
\]
while for $x=y$ we have $H^\bullet(\Ss^d-\{x\}) = \k$ and $H^\bullet_c(\Ss^d-\{x\}) = \k[-d]$.

Further, the cotwist $T_0$ is the shift by $(-d)$ tensored with $H^d(\Ss^d, \k)$, 
 the
1-dimensional vector space of orientations of $\Ss^d$. 
  So it is an equivalence and (SF2) is also satisfied.
 \qed

\begin{ex}\label{ex:p1-spec}
Note the particular case $d=2$, when $\Ss^2=\CC\PP^1$. Proposition \ref{prop:p1-sph} implies that for any
$\PP^1$-fibration $q: Z\to Y$ of complex algebraic varieties, $q^{-1}: D^b(Y)\to D^b(Z)$ is a spherical functor.
Another class of examples is provided by quaternionic geometry, since $\HH\PP^1 = \Ss^4$. 
\end{ex}

Let us now mention some ``coherent''  examples.

\begin{ex}
Recall that an $n$-dimensional smooth projective variety $Z$ over $\k$ is called {\em Calabi-Yau}  (in the strict sense), if
\[
H^i(Z, \Oc_Z)=\begin{cases}
\k, & \text{ if } i=0, \text{ or } i=n,
\\
0, & \text{ otherwise}. 
\end{cases}
\]
Let $X$ be a smooth algebraic variety over $\k$, and $q: Z\to X$ be a smooth proper family of Calabi-Yau
manifolds. Then the pullback functor $q^*: D^b_\coh(X) \to D^b_\coh(Z)$ is a spherical functor. 
The proof is similar to that of Proposition \ref{prop:p1-sph}. 
\end{ex}

\begin{ex}\label{ex:YDX}
(\cite{A1})
Consider a diagram 
\[
Y\overset{\rho}\lla D\overset{i}\hra X
\]
of smooth  complex varieties $X, D, Y$, where $\rho$ is a $\PP^1$-bundle and 
$i$ an embedding of a divisor. We then have a diagram of adjoint functors 
\[
\begin{matrix}\ &\overset{L = \rho_! i^*}\lla &\ \\
\Dc_0:=  D^b_\coh(Y)& \overset{S = i_* \rho^*} \lra & D^b_\coh(X) \\
& \overset {R=\rho_* i^!} \lla &
\end{matrix}
\]
and $R = T_0L$ where 
$T_0 = \Cone\{\Id_{\Dc_0}\Rightarrow RS\}$. 
\end{ex}
 
\begin{lem}\label{lem:A-div}
(R. Anno, \cite{A1}). $S$ is spherical iff the intersection 
index of $D$ with a generic fiber of $\rho$ is $(-2)$.  \qed
\end{lem}

\begin{ex}\label{ex:p1-coh}
A particular case of Example \ref{ex:YDX} and Lemma \ref{lem:A-div}
 is obtained for 
\[
Y=\pt \buildrel \rho\over\lla D=\PP^1 \buildrel i\over\hra X=T^*\PP^1
\]
The corresponding spherical functor can be seen as a ``quasi-classical approximation'' to that in
Example \ref{ex:p1-spec}. 
\end{ex}

\begin{rem}

More precisely, for a complex 
 algebraic manifold $Z$, the coherent derived category $D^b_\coh(T^*Z)$ can be thought of as a
 quasi-classical approximation to the derived category $D^b(Z)$ of arbitrary sheaves on $Z$.
Indeed, passing to solutions of $\Dc$-modules gives a functor between derived categories in the first
line of the following table (a functor that restricts to the Riemann-Hilbert equivalence between constructible and
holonomic regular derived categories). It can be compared with the second line which is an instance of
Serre's theorem for the affine morphism $p: T^*Z\to Z$. 
    
  \vskip .3cm

  \begin{tabular}{ p{5cm} | p{5cm}  }
  (Arbitrary) sheaves on $Z$ & Coherent  $\Dc_Z$-modules on $Z$
  \\
  \hline 
  Coherent sheaves on $T^*Z$ &  Coherent  $p_*\Oc_{T^*Z} = \on{gr}(\Dc_Z)$-modules   on $Z$  \\ 
  \end{tabular}

\vskip .3cm

Note that the functor $S=q^{-1}$  on sheaf-theoretic derived categories in Example \ref{ex:p1-spec} matches, after being interpreted
in terms of $\Dc$-modules and passing to the associated graded, the functor $S=i_*\rho^*$ on coherent derived categories in
Example \ref{ex:p1-coh}. 

\end{rem}

 \vfill\eject

\section {\bf Disk with several marked points}\label{sec:disk-several}

\noindent {\bf A. Quiver description of perverse sheaves.} 
Let $B=\{b_1, \cdots, b_n\}$ be a finite set of marked points in the unit disk $\Delta$ and $\Perv(\Delta, B)$
be the category of perverse sheaves on $\Delta$ with possible singularities at $B$. We then have the following
\cite[Prop. 1.2]{gelfand-MV}.

\begin{prop}\label{prop:GMV}
$\Perv(\Delta, B)$ is equivalent to the category $\Pc_n$    of diagrams formed by vector spaces $\Psi, \Phi_1, \cdots, \Phi_n$
and linear maps 
     $\xymatrix{
 \Phi_i \ar@<.4ex>[r]^{v_i}&\Psi \ar@<.4ex>[l]^{u_i}
 }$
    such that each 
  $T_{\Psi, i}: = \Id_\Psi - v_iu_i$ is an isomorphism. 
 \end{prop}
 
 The category $\Pc_n$ can be seen as an amalgamation of $n$ copies of the category $\Pc_1$ from Theorem \ref{thm:ggm}. 
 To construct an equivalence, we choose a base point $b\in\partial\Delta$ and a ``system of cuts" $K$,
 i.e., a set of simple arcs $\{K_1, \cdots, K_n\}$ with $K_i$ connecting $b$ with $b_i$ and
 with different $K_i$ meeting only at 
 a small  common interval near $b$,
 see Fig. \ref{fig:disk-2}. Notationally, we view $K$ as the union  $K = \bigcup K_i\subset\Delta$.

   \bef
\centering
\btp[scale=.4, baseline=(current  bounding  box.center)]
 
 \node (1)  at (-1,2.5){}; 
\fill (1) circle (0.15);   

\node (2)  at (-2.5,1){};
\fill (2) circle (0.15);  

\node (n)  at (-2,-3){}; 
\fill (n) circle (0.15);

  \draw (0,0) circle (7cm);

   \node (p) at (6,0){}; 
\fill (p) circle (0.15);

 \node (b) at (7,0){}; 
\fill (b) circle (0.15);

\draw[line width = .3mm] (6,0) .. controls (2,1) .. (-1,2.5); 

\draw[line width = .3mm] (6,0) .. controls (0,-1) .. (-2.5,1); 

\draw[line width = .3mm] (6,0) .. controls (2,6) and (-9,6)   .. (-2,-3); 

\draw[line width = .3mm] (6,0) -- (7,0);

\node at (7.7, .5) {$b$}; 
 
\node at (6,6) {$\Delta$}; 

\node at (-1.6, 2.8) {$b_1$}; 

\node at (-3.3, 1.5){$b_2$}; 

\node at (-2.2, -1){$\vdots$}; 

\node at (-1, -3){$b_n$}; 

\node at (-4,-2){$K_n$}; 

\node at (2,2){$K_1$}; 

\node at (0.5,0){$K_2$};

\etp

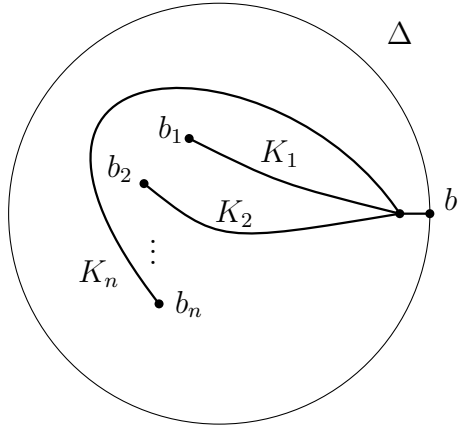
\caption{ Equivalence depending on a system of cuts.  }
\label{fig:disk-2}
\enf
Given $K$, we have an equivalence $F_K: \Perv(\Delta, B)\to \Pc_n$ sending $\Fc$ to 
\be\label{eq:F-K}
\Phi_i^K (\Fc) = \ul\HH^1_K(\Fc)_{b_i}, \quad \Psi^K(\Fc) = \Fc_b \simeq \ul\HH^1_{K_i}(\Fc)_{b} = \ul\HH^1_{K}(\Fc)_b. 
\ee

\begin{rem}\label{rem:patt-2}
 We have the following elementary statements, continuing the pattern of Remark \ref{rem:patt-1}.
 \begin{enumerate}
\item[(1)]  $\ul\HH^i_K(\Fc) = 0$  { for }  $i\neq 1$.

\item[(2)] Each   stalk of $\Rc(\Fc)$, considered as a functor $\Perv(\Delta,B)\to\Vect_\k$, is an exact functor which takes
Verdier duality to vector space duality.
\end{enumerate}
\end{rem}

\begin{rem}
The  {\em space}  of cohomology with support $\HH^1_K(\Delta, \Fc)$ can be seen as ``uniting'' all the spaces of vanishing cyclies
$\Phi_i^K(\Fc)$. As common in singularity theory, we can
 imagine  that $\Fc$ obtained as a deformation of a perverse sheaf  $\Gc\in\Perv(\Delta, 0)$
  with only one (but more complicated)  singular point at $0$.
 Then $\HH^1_K(\Delta, \Fc)$ recovers $\Phi(\Fc)$. Note that a categorification of the space of vanishing cycles
 of an isolated singular point  of a function
  is provided by the Fukaya-Seidel category \cite{FSbook}, and the method of construction adopted in {\em loc. cit.} 
  uses precisely a deformation into several Morse critical points.  Therefore,  the space $\HH^1_K(\Delta, \Fc)$ for a
  perverse sheaf $\Fc\in\Perv(\Delta, B)$,   can be seen as a de-categorified analog of the Fukaya-Seidel category. 
\end{rem}

\vskip .2cm

\noindent 
The equivalence $F_K$  depends only on the isotopy class of $K$.
Unlike the one point case, there are now many such classes, forming a set which we denote $\Cc$.  
It
is acted upon simply transitively by the Artin braid group 
\[
\Br_n \,\,=\,\, \big\langle s_1, \cdots, s_{n-1} \bigl| \,\,\, s_i s_{i+1} s_i = s_{i+1} s_i s_{i+1}\big\rangle.
\]
Indeed,
\[
\Br_n \,\,=\,\,\pi_0 \on{Diff}^+(\Delta; B, b)
\]
is the group of isotopy classes of  diffeomorphisms of $\Delta$, preserving orientation, preserving $b$ as a point and $B$ as a set.
The equivalences $F_K$ for different $K\in\Cc$ are connected by self-equivalences $f_\sigma$  of $\Pc_n$:
\be\label{eq:self}
\xymatrix{
\Pc_n \ar[rr]^{f_\sigma} && \Pc_n, 
\\
& \Perv(\Delta, B)
\ar[ul]^{F_K} \ar[ur]_{F_{\sigma(K)}}&
}
\quad \sigma\in\Br_n.
\ee
The self-equivalence $f_{s_i}$ corresponding to a generator $s_i$ of $\Br_n$,  is given by  
\cite[Prop. 1.3]{gelfand-MV}:
\be\label{eq:PL0}
\begin{gathered}
f_{s_i}\bigl( \Psi, \Phi_j, u_i, v_j\bigr) \,=\, \bigl( \Psi, \Phi'_j, u'_j, v'_j\bigr),\\
\Psi'=\Psi, \,\,\, \Phi'_j = \Phi_j, \,\,u'_j=u_j,\,\, v'_j=v_j, \quad j\neq i, i+1, \\
\Phi'_{i+1}=\Phi_i, \,\,\Phi'_i = \Phi_{i+1}, \\
u'_i = u_{i+1}, \,\, v'_i = v_{i+1}, \,\,\, u'_{i+1} = u_i T_{\Psi, i+1}, \,\, v'_{i+1} = T_{\Psi, i+1}^{-1} v_i. 
\end{gathered}
\ee

\begin{rems}\label{rems:coord}
 (a) Note that  perverse sheaves being a topological concept, the group
  $\on{Diff}^+(\Delta; B, b)$ naturally acts on the category $\Perv(\Delta, B)$ from the first principles, 
 the action descending to that of $\Br_n$. 
 
 \vskip .2cm

 (b) We can  turn Proposition \ref{prop:GMV}  and formulas \eqref{eq:PL0} 
around to produce an intrinsic (i.e., not tied to any particular $K$ and manifestly $\Br_n$-equivariant) {\em definition}
 of $\Perv(\Delta, B)$  which does not appeal to any pre-existing concept of a perverse sheaf. 
 More precisely, we can {\em define} an object  $P\in\Perv(\Delta, B)$ to be a system of objects $P_K\in\Pc_n, K\in \Cc$ and compatible
 isomorphisms $f_\sigma(P_K) \to P_{\sigma(K)}$, $\sigma\in \Br_n$ so that each particular $P_K$ is
 just  a  particular ``shadow'' of a more intrinsic object $P$.

\end{rems}

\vskip .3cm

\noindent {\bf B. Schobers on a disk with several marked points.} 
To give an ``invariant"  definition of a perverse Schober on $(\Delta, B)$, we adopt the approach of
Remark \ref{rems:coord}(b). That is, for each system of cuts $K\in\Cc$ we define a {\em $K$-coordinatized Schober}
to be a system of $n$ spherical functors with a common target
\[
\Sen_K\,\,=\bigl\{  S_i: \Dc_i \lra \Dc, \quad i=1, \cdots, n\bigr\}. 
\]
Each  $\Sen_K$ gives rise to spherical reflection functors $T_i = \Cone\{S_iR_i\Rightarrow \Id_\Dc\}$, $i=1, \cdots, n$. 
According to our convention on working with triangulated categories in terms of dg-enhancements, see Appendix,
all $K$-coordinatized Schobers form an $\infty$-category which we denote $\on{Sch}_K(\Delta, B)$. We  define
equivalences $\fb_\sigma: \Sch_K(\Delta, B)\to \Sch_{\sigma(K)}(\Delta, B)$, $\sigma\in\Br_n$ on generators $s_i\in\Br_n$ by
the direct analog of \eqref{eq:PL0}: 

\be\label{eq:PL1}
\begin{gathered}
{\bf f}_{s_i} \bigl\{  S_i: \Dc_i \lra \Dc\bigr\}\,\,=\,\, \bigl\{  S'_i: \Dc'_i \lra \Dc' \bigr\}, 
 \\
\Dc'=\Dc, \,\,\, \Dc'_j = \Dc_j, \,\,S'_j=S_j,\,\,   \quad j\neq i, i+1, \\
\Dc'_{i+1}=\Dc_i, \,\,\Dc'_i = \Dc_{i+1}, \\
  S'_i = S_{i+1}, \,\,\,  S'_{i+1} = T_{ i+1}^{-1} S_i.
\end{gathered}
\ee
We then extend to arbitrary $\sigma\in\Br_n$ by verifying the braid relations for the $\fb_i$ which is done in exactly the
same way as for \eqref{eq:PL0}.

By definition, a {\em perverse Schober on $(\Delta, B)$} is a system $\Sen = (\Sen_K)_{K\in\Cc}$ of coordinatized Schobers
and compatible identifications $\fb_\sigma(\Sen_K) \to\Sen_{\sigma(K)}$. The datum $\Sen_K$ will be referred to as the
{\em $K$-shadow} of $\Sen$.  We denote by $\Sch(\Delta, B)$ the $\infty$-category of perverse Schobers on $(\Delta, B)$.

\vskip .3cm

\noindent {\bf C. The Picard-Lefschetz formula.} Underlying classical Picard-Lefschetz theory,  there is a general statement about perverse
sheaves on a disk which we now formulate.

Let $\Fc\in\Perv(\Delta, B)$. In the approach of \eqref{eq:F-K}, ``the''   space  of vanishing cycles of $\Fc$ at  some $b_i\in B$ can be defined in terms
of a small segment of an arc terminating in $b_i$ (which does not, {\em a priori}, have to be a part of a system of cuts). 
Let now $\gamma$ be a simple arc joining two marked points $b_i$ and $b_k$ and not passing through any other marked points,
as in Figure \ref{fig:PL}.

 \bef
\centering
\btp[scale=.4, baseline=(current  bounding  box.center)]

 \draw (0,0) circle (8cm); 
   
   \node (p) at (7,0){};

   \node (b) at (8,0){}; 
\fill (b) circle (0.15);

\node at (8.7, .5) {$b$}; 
 
\node at (7,7) {$\Delta$}; 

\node (i) at (-3,3.5){};
\fill(i) circle (0.15); 

\node (k) at (-3.5,-3){};
\fill (k) circle (0.15); 

\node (j) at (-1,.7){};
\fill (j) circle (0.15); 

\draw [dotted, line width = .3mm] (7,0) -- (8,0); 

\draw [line width = .3mm]  (-3,3.5) .. controls (3,1) .. (-3.5,-3); 

\draw [line width = .3mm]  (-3,3.5) .. controls (-4,1) .. (-3.5,-3); 

\draw [line width = .3mm]  (-3,3.5) .. controls (-2,2) .. (-1,0.7); 

\draw [line width = .3mm]  (-3.5,-3) .. controls (-3,-2) .. (-1,0.7); 

\draw [dotted, line width = .3mm]  (-3.5,-3) .. controls (-1,-2) .. (7,0); 

\draw [dotted, line width = .3mm]  (-1, 0.7) .. controls (0,0) .. (7,0); 

\draw [dotted, line width = .3mm]  (-3, 3.5) .. controls (2,3) .. (7,0); 

\node at (-3.5, 4.5){$b_i$};

\node at (-3.5, -4){$b_k$};

\node at (-2, 0.7){$b_j$};

\node at (-4.6,0){$\gamma$};

\node at (-2.7,-.8){$\alpha$};

\node at (-2.6,2){$\beta$};

\node at (1.9,1.4){$\gamma'$};

\node at (3.2,3){$K$}; 
\etp

\caption{ The Picard-Lefschetz situation.  }
\label{fig:PL}
\enf

We can then  define the spaces of vanishing and nearby cycles of $\Fc$ relative to $\gamma$:
\[
\begin{gathered}
\Phi_{i,\gamma} = \ul{\HH}^1_\gamma(\Fc)_{b_i}, \quad \Phi_{k, \gamma} = \ul{\HH}^1_\gamma(\Fc)_{b_k},\\
\Psi_\gamma =  \ul{\HH}^1_\gamma(\Fc)_{\on{gen}}\,\text{ (generic stalk)}. 
\end{gathered}
\]
which are connected by the maps
\[
  \xymatrix{
 \Phi_{i,\gamma}  \ar@<.4ex>[r]^{v_{i,\gamma}}&\Psi_\gamma  \ar@<.4ex>[l]^{u_{i,\gamma}} 
  \ar@<-.4ex>[r]_{u_{k,\gamma}}
 &
 \Phi_{k, \gamma} 
  \ar@<-.4ex>[l]_{v_{k,\gamma}} .
 }
\]
The definition of these maps is similar to  \S \ref{sec:disk}A: the maps $v$ are generalization maps, and the maps
$u$ are obtained by counterclockwise continuation, as in \eqref{eq:map-u}. 
We define the {\em transition map} along $\gamma$ as
\[
M_{ik}(\gamma) \,\,=\,\,  u_{k,\gamma} \circ v_{i,\gamma}: \Phi_{i,\gamma}\lra\Phi_{k, \gamma}.
\]
The data of the $\Phi_i(\Fc)$, $i=1, \cdots, n$ as local systems on the circles around $b_i$ and
of all the $M_{ik}(\gamma)$ describe the image of $\Fc$ in the localization of $\Perv(\Delta, B)$ by the subcategory
of constant sheaves \cite[\S 2]{gelfand-MV}.  

\vskip .2cm

One can say that ``abstract Picard-Lefschetz theory'' is the study of
 how $M_{ik}(\gamma)$ changes when we replace $\gamma$ by a different
(non-isotopic) arc $\gamma'$. 

More precisely,  assume that $\gamma'$ is obtained from $\gamma$ by an ``elementary move'' past another marked point $b_j$ so that
the bigon formed by $\gamma'$ and $\gamma$ contains $b_j$ and two arcs $\alpha$ and $\beta$, but no other
marked points  as in Fig. \ref{fig:PL}. 
In particular, we assume that the closed path obtained by following $\gamma$ from $b_i$ to $b_k$ and then $\gamma'$ from $b_k$ to $b_i$,
has orientation compatible with the standard (counterclockwise) orientation of $\Delta$. 
Note that  homotopy inside the bigon  and the {\em clockwise} rotation around $b_j$ give identifications 
\be\label{eq:PL-identif}
\Phi_{i,\gamma} \simeq \Phi_{i,\beta} \simeq \Phi_{i, \gamma'}, \quad \Phi_{k,\gamma}\simeq\Phi_{k, \alpha}\simeq \Phi_{k,\gamma'},
\quad \Phi_{j,\beta} \simeq \Phi_{j, \alpha}, 
\ee
so we can consider them as single spaces denoted by $\Phi_i, \Phi_k$ and $\Phi_j$ respectively. 

\begin{prop}[Picard-Lefschetz formula for perverse sheaves]\label{prop:PL-form}
We have the equality of linear operators $\Phi_i\to\Phi_k$: 
\[
M_{ik}(\gamma') = M_{ik}(\gamma) + M_{jk}(\alpha) M_{ij}(\beta). 
\]

\end{prop}

This statement  is a version of  \cite[Prop. 2.4]{gelfand-MV}, formulated in a more invariant way and
without localizing by constant sheaves.
It holds  for perverse sheaves on any oriented surface. It is convenient to give two proofs of Proposition \ref{prop:PL-form}. 

\vskip .2cm

\noindent {\sl Invariant proof:}  
To eliminate the need for the first two identifications in \eqref{eq:PL-identif}, let us deform the paths
$\alpha,\beta,\gamma, \gamma'$ so that:
\begin{itemize}
\item $\gamma, \beta, \gamma'$ have a common segment $[b_i, b'_i]$ near $b_i$.

\item $\gamma, \alpha, \gamma'$ have a common segment $[b'_k, b_k]$ near $b_k$. See Fig. \ref{fig:PL-def}. 
\end{itemize}
We denote by $\ol\gamma, \ol\gamma'$ the parts of $\gamma$ and $\gamma'$ lying between $b'_i$ and $b'_k$,
by $\ol\beta$ the part of $\beta$ between $b'_i$ and $b_j$, and by $\ol\alpha$ the part of $\alpha$  between $b_k$
and $b'_j$. 

 \bef
\centering
\btp[scale=.4, baseline=(current  bounding  box.center)]
 
 \node (b') at (0,0){};
 \fill (b') circle (0.15); 
 
 \node (b) at (-3,4){};
 \fill (b) circle (0.15); 
 
 \draw [line width = .3mm] (0,0) -- (-3,4); 
 
 \draw  (0,0) -- (-3,-3); 
 \draw[dotted, line width = .3mm] (-3,-3) -- (-4,-4); 
 
 \draw (0,0) -- (0,-3); 
 
 \draw[dotted, line width = .3mm] (0,-3) -- (0,-4);

 \draw (0,0) -- (4,-3); 
 
 \draw[dotted, line width = .3mm] (4,-3) -- (5,-3.75);  
 
 \node at (-2, 4.5){$b_i$}; 
 
 \node at (1, .5) {$b'_i$}; 
 
 \node at (-4, -3) {$\ol\gamma$}; 
 
 \node at (-1, -3) {$\ol\beta$}; 
 
 \node at (4, -2) {$\ol\gamma'$}; 
 \etp
 \quad\quad\quad\quad
\btp[scale=.4, baseline=(current  bounding  box.center)]

 \node (b') at (0,0){};
 \fill (b') circle (0.15); 
 
  \node (b) at (2,-4){};
 \fill (b) circle (0.15); 
 
  \draw [line width = .3mm] (0,0) -- (2,-4); 
  
  \draw (0,0) -- (-3,3); 
  \draw[dotted, line width = .3mm] (-3,3) -- (-4,4);  
  
  \draw (0,0) -- (1,3); 
   \draw[dotted, line width = .3mm] (1,3) -- (1.33,4);

  \draw (0,0) -- (4,2); 
 \draw[dotted, line width = .3mm] (4,2) -- (5, 2.5);  
 
 \node at (-1,-.5) {$b'_k$}; 
 
 \node at (3, -4.5) {$b_k$}; 
 
 \node at (-3.5, 2) {$\ol\gamma$}; 
 
 \node at (0,3) {$\ol\alpha$}; 
 
 \node at (4,1) {$\ol\gamma'$}; 

\etp

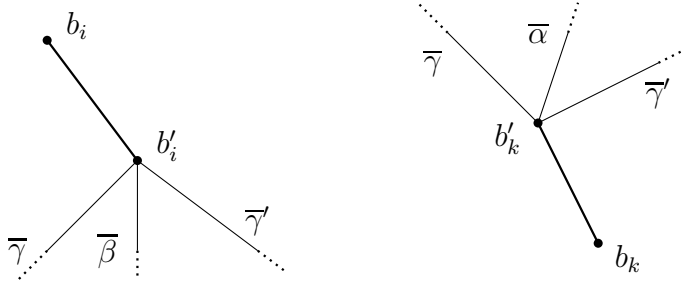
\caption{ The Picard-Lefschetz situation, deformed.  }
\label{fig:PL-def}
\enf

The points $b'_i, b'_k$ being smooth for $\Fc$, both sides of our putative equality factor through the maps
\[
v_{[b_i, b'_i]}: \Phi_i \lra \Fc_{b'_i}, \quad u_{[b'_k, b_k]}: \Fc_{b'_i} \lra \Phi_k. 
\]
Our statement would therefore follow from the next lemma which is an invariant version of the statement that
the map \eqref{eq:t-psi} is indeed the monodromy around 0.

 \begin{lem}
 We have an equality of operators $\Fc_{b'_i}\to \Fc_{b'_k}$:
 \[
 T_{\ol\gamma'} \,\,=\,\, T_{\ol\gamma} + v_{\ol\alpha} R u_{\ol\beta},
 \]
 where:
 \begin{itemize}
 \item[(1)]  $u_{\ol\beta}: \Fc_{b'_i}\to \ul\HH^1_{\ol\beta}(\Fc)_{b_j}$ is the coboundary of the counterclockwise continuation map,
 cf. \eqref{eq:map-u}. 
 
 \item[(2)] $R: \ul\HH^1_{\ol\beta}(\Fc)_{b_j} \to \ul\HH^1_{\ol\alpha}(\Fc)_{b_j}$ is the identification obtained by deforming, by
 clockwise rotation aroung $b_j$, the path $\ol\beta$ into the path $\ol\alpha$, i.e., $R$ is the third identification in \eqref {eq:PL-identif}.
 
 \item[(3)] $v_{\ol\alpha}: \ul\HH^1_{\ol\alpha}(\Fc)_{b_j}\to \ul\HH^1_{\ol\alpha}(\Fc)_{\on{gen}} \simeq\Fc_{b'_k}$
 is the generalization map of the sheaf $\ul\HH^1_{\ol\alpha}(\Fc)$. 
 
  \end{itemize}
 
 \end{lem}
 
 \noindent {\sl Proof of the lemma:} Let $U$ be a small disk around $b_j$. We fix $a\in\Fc_{b'}$ 
 and compare the two sides of the putative equality when applied to $a$. For this, 
 let $s_a\in\Gamma(U-\ol\beta, \Fc)$
 be the section obtained by continuing $a$ on the left side of $\ol\beta$ towards and around $b_j$. Then $u_{\ol\beta}(a)$
 is equal to the class of $s_a$ in $\ul\HH^1_{\ol\beta}(\Fc)_{b_j}$. Next, $Ru_{\ol\beta}(a)$ is similarly  represented by the section
 $t_a$ obtained from $s_a$ by continuously moving the branch cut clockwise from $\ol\alpha$ to $\ol\beta$. This means that $t_a=s_a$
 on the left side of $\ol\alpha \cup\ol\beta$. Finally, $v_{\ol\alpha} R u_{\ol\beta}(a)$ is obtained as the difference of the two boundary
 values if $t_a$ when continued along both sides of $\ol\alpha$ all the w ay to $b'_k$. It remains to notice that these
 boundary values are equal, in virtue of the above, to $T_{\ol\gamma}(a)$ and $T_{\ol\gamma'}(a)$. \qed.

 \vskip .3cm

 \noindent {\sl Proof using shadows:}   Choose a system of cuts $K$ (depicted by dotted lines in Fig. \ref{fig:PL})
 adopted to our situation.  We assume 
   that the arcs $K_i, K_j$ and $K_k$ are
positioned as in the figure, i.e., that $\gamma$ together with $K_i$ and $K_k$ form a triangle containing $K_j$, $\alpha$ and $\beta$
and not containing any other marked points. We  orient each  $K_\nu$ to run from $b$ to $b_\nu$. 
Consider   the quiver $F_K(\Fc) = (\Psi, \Phi_i, u_{i}, v_{i})$.   
Note that we have isotopies of oriented paths rel. $B$:
\[
\begin{gathered}
\alpha \sim K_k*K_j^{-1}, \,\,\,\beta \sim K_j*K_i^{-1}, \,\, \gamma'\sim K_k*K_i^{-1}, \\
\gamma\sim K_k * \partial\Delta * K_i^{-1}.
\end{gathered}
\]
Here  $*$ means composition of the paths and $\partial\Delta$ is the boundary circle of $\Delta$, oriented anticlockwise and
run from $b$ to $b$. We can use these isotopies to calculate the transition maps, obtaining
 \[
\begin{gathered}
M_{ij}(\beta) =  u_{j,K}v_{i,K}, \quad M_{jk}(\alpha) =  u_{k,K}v_{j,K}, \quad
M_{ik}(\gamma') =  u_{k, K} v_{i,K}, \\
 \quad M_{ik}(\gamma) = u_{k,K} T_{j, \Psi} v_{i,k}, 
\end{gathered}
\]
and the claim follows from the identity $T_{j,\Psi}= \Id-v_i u_j$. \qed

\vskip .3cm

\noindent {\bf D. The Picard-Lefschetz triangle.}
Let $\Sen$ be a perverse Schober on $(\Delta, B)$. The transition maps constructed in n$^\circ$C categorify to functors between
triangulated categories. More precisely, if $\gamma$ is an oriented arc joining $b_i$ and $b_k$ as above, then we have
a diagram of triangulated categories and spherical functors
which depends ``canonically" (i.e., up to a contractible set of choices) only on the isotopy class of $\gamma$ rel. $B$:
\be\label{eq:maps-gamma-cat}
  \xymatrix{
 \Phi_{i,\gamma}(\Sen)  \ar@<.0ex>[r]^{S_{i,\gamma}}&\Psi_\gamma(\Sen)   
   &
 \Phi_{k, \gamma} (\Sen) 
  \ar@<-.0ex>[l]_{S_{k,\gamma}} .
 }
\ee
To define it, we choose $K'\in\Cc$ so that $K'_i$ (oriented from $b$ to $b_i$), together with $\gamma$ and $(K'_k)^{-1}$ form
a positively oriented triangle not containing any other marked points. (Note: the choice of $K$ as in Fig.\ref{fig:PL}
is not good.)
We consider the $K'$-shadow $\Sen_{K'}= \bigl\{S_\nu^{K'}: \Dc_\nu^{K'}\to\Dc^{K'}\bigr\}_{\nu=1}^n$ of $\Sen$ 
and define
\[
\begin{gathered}
\Psi_\gamma(\Sen) = \Dc^{K'}, \,\,\,\Phi_{i,\gamma}(\Sen) = \Dc_i^{K'}, \,\,\, \Phi_{k,\gamma}(\Sen) = \Dc_k^{K'}, 
\\
S_{i, \gamma} = S_i^{K'}, \,\,\, S_{k, \gamma} = S_k^{K'}. 
\end{gathered}
\]
It is straightforward to see (by looking at the subgroup in $\Br_n$ permuting all $K'$ with our property) that
this definition is indeed ``canonical'' in the sense  described. 
So we consider the data \eqref{eq:maps-gamma-cat} as intrinsically associated to $\Sen$ and $\gamma$.

Denote $R_{k,\gamma}$   the right  adjoint to the spherical functor $S_{k,\gamma}$.
  Define now the {\em transition functor}
\[
\Mc_{ik}(\gamma) = R_{k,\gamma} \circ S_{i,\gamma}: \, \Phi_{i, \gamma}(\Sen) \lra\Phi_{k,\gamma}(\Sen). 
\]
Consider now a situation when we have arcs $\gamma, \gamma',\alpha, \beta$
 depicted in Fig. \ref{fig:PL}. Using identifications of  categories similar to
\eqref{eq:PL-identif}, we 
  can speak about
triangulated categories $\Phi_i(\Sen)$ and $\Phi_k(\Sen)$.

\begin{prop}[Picard-Lefschetz triangle]  We have a canonical triangle of  exact functors $\Phi_i(\Sen)\to \Phi_k(\Sen)$: 
\[
\Mc_{ik}(\gamma) \lra \Mc_{ik}(\gamma') \lra \Mc_{jk}(\alpha) \circ\Mc_{ij}(\beta) \lra\Mc_{ik}(\gamma)[1]. 
\]
\end{prop}

\noindent {\sl Proof:} This is  obtained identically to the ``shadow" proof of Proposition \ref{prop:PL-form} with the identity
$T_{j,\Psi} = \Id-v_ju_j$ replaced by the triangle coming from the fact that $T_j = \Cone\{ \Id \to R_jS_j\}$. \qed

 \vskip .3cm
 
 \noindent{\bf E.  Schobers on a Riemann surface.} Let $\Sigma$ be an oriented topological surface, possibly with
 boundary $\partial\Sigma$, and $B\subset \Sigma$ is a finite set not meeting $\partial\Sigma$. One can then
 define a perverse Schober on $(\Sigma, B)$ by decomposing $\Sigma$ as $\Delta\cup_C  U$,
 where $\Delta\subset\Sigma$ is a closed disk with boundary circle $C=\partial\Delta$ which contains
 all points of $B$, and $U$ is the closure of $\Sigma-\Delta$. Then a perverse Schober $\Sen$
 is, by definition,  a datum of:
 \begin{itemize}
  \item[(1)] A perverse Schober $\Sen_\Delta$ on $(\Delta, B)$, defined as in n$^\circ$B.
  
  \item[(2)] A local system of triangulated categories $\Sen_U$ on $U$ identified with $\Sen_\Delta$ over $C$. 
 
 \end{itemize}
While one can work with objects thus defined (for instance, one can 
construct  transition functors and Picard-Lefschet triangles for arcs not necessarily
contained in $\Delta$), a more intrinsic definition is desirable.

\vfill\eject

\section{Spherical functors and spherical pairs}

\noindent{\bf A. Symmetric description of $\Perv(\Delta, 0)$.}
 In \cite[\S 9]{KS} we have given a different quiver description of $\Perv(\Delta,0)$
 obtained as a particular case of a general result for real hyperplane arrangements:

\begin{prop}\label{prop:C-real}
 The category $\Perv(\Delta,0)$ is equivalent to the category formed by diagrams of vector spaces
  \be\label{eq:diag-real}
 \xymatrix{
E_-\ar@<-.7ex>[r]_{\delta_-}& E_0
\ar@<-.7ex>[l]_{\gamma_-}
\ar@<.7ex>[r]^{\gamma_+}& E_+
\ar@<.7ex>[l]^{\delta_+}
}
\ee
satisfying the two following conditions:
\begin{enumerate}
\item[(1)]  $\gamma_-\delta_- = \on{Id}_{E_-}, \,\,\,\gamma_+
\delta_+ = \on{Id}_{E_+}$.
\item[(2)] The maps $\gamma_-\delta_+: E_+\to E_-, \,\,\, \gamma_+\delta_-: E_-\to E_+ $
are invertible. 
\end{enumerate}
 \end{prop}
 This is obtained by choosing not one but two base points $b_+, b_-\in\partial\Delta$ and considering  a cut $K$
 which joins $b_+$ with $b_-$ and passes through $0$, as in  depicted in Fig. \ref{fig:disk-4}.

    \bef
\centering
\btp[scale=.4, baseline=(current  bounding  box.center)]
\node (0) at (0,0){}; 
\fill (0) circle (0.15);
   
 \draw (0,0) circle (5cm);

   \node (b+) at (5,0){}; 
\fill (b+) circle (0.15);

 \node (b-) at (-5,0){}; 
\fill (b-) circle (0.15);

\draw[line width = .3mm] (5,0) ..  controls (1,1)  and (-1, -1)  .. (-5,0); 

\node at (-.7, .7) {$0$}; 
\node at (5.7, .5) {$b_+$}; 
\node at (1.7,1) {$K$}; 
\node at (4.5,4.5) {$\Delta$}; 
\node at (-5.7, .5) {$b_-$};

\etp

\caption{ A symmetric cut.  }
\label{fig:disk-4}
\enf
\noindent The spaces $E_\pm, E_0$ are obtained as  the stalks of  the sheaf $\ul\HH^1_K(\Fc)$ at $b_\pm$ and $0$ respectively,
the maps $\gamma_\pm$ are the generalization maps for this sheaf, and the $\delta_\pm$ can be obtained by duality.

\begin{rems}\label{rem:patt-3}
(a)  Note that in this description the maps $P_+=\delta_-\gamma_-$ and $P_+ = \delta_+\gamma_+$
are projectors in $E_0$, that is $P_\pm^2=P_\pm$. We can consider $E_\pm$ as  subspaces
in $E_0$ which are the images of $P_\pm$. 

\vskip .2cm

(b) The pattern of Remarks \ref{rem:patt-1} and  \ref{rem:patt-2} continues here: $\ul\HH^i_K(\Fc)=0$ for $i\neq 1$ and
each stalk of $\ul\HH^1_K$ considered as a functor into $\Vect_\k$, takes Verdier duality into vector space duality. 

\end{rems}

 \vskip .3cm
 
 \noindent {\bf B. Reminder on semi-orthogonal decompositions.} 
 A categorification of diagrams \eqref{eq:diag-real}
  is naturally formulated in the language of semi-orthogonal decompositions
 of triangulated categories \cite{BK} \cite{kuznetsov} which we now recall.
 
Let  
$\Bc$ be a full triangulated subcategory of a triangulated category $\Ac$,
 with $i=i_\Bc: \Bc\to\Ac$ being the embedding functor.  We denote by
 \[
\begin{gathered}
^\perp \Bc =   \{A\in\Ac| \Hom(A,B)=0 \,\,\forall B\in\Bc\}, \\
 \Bc^\perp  = \{A\in\Ac| \Hom(B, A)=0 \,\,\forall B\in\Bc\}
 \end{gathered}
\]
the left and right orthogonals to $\Bc$
  We say that $\Bc$ is {\em left admissible},
 resp. {\em right admissible}, if $i$ has a left adjoint $^*i$, resp. a right adjoint $i^*$. 
 If $\Bc$ is left (resp. right)
 admissible, we have a semi-orthogonal decomposition  
 \[
\Ac = \langle ^\perp\Bc, \Bc\rangle , \quad \text{ resp.} \,\,\,
 \Ac = \langle \Bc, \Bc^\perp\rangle
\]
which means  
 that each object $A\in \Ac$ is included in functorial exact triangles
\[
\begin{gathered}
\hskip 1cm C\lra A\lra B\to  C[1], \quad C\in ^\perp\Bc, \,\,\,\, B = {} ^*i(A) \in \Bc,\\
 \text { resp.}\,\,\,\,  B'\lra A\lra D'\to B'[1], \quad B'= i^*(A) \in \Bc, \,\,\,\,  D'\in{}  \Bc^\perp.  
\end{gathered}
\]
In particular, $^\perp\Bc = \Ker({} ^*i)$, resp. $\Bc^\perp = \Ker(i^*)$. 
We will call $^*i$ the {\em projection on $\Bc$ along $^\perp\Bc$}  and
$i^*$ the {\em projection on $\Bc$ along $\Bc^\perp$}. 

If $\Bc$ is left admissible, then $^\perp\Bc$ is right admissible, and $(^\perp \Bc)^\perp = \Bc$.
Similarly, for a right admissible $\Bc$ we have that $\Bc^\perp$ is left admissible and $^\perp(\Bc^\perp)=\Bc$. 
 We call $\Bc$ {\em admissible} if it is both left and right admissible. 
 For an admissible $\Bc$ we have an equivalence
 \[
 M_\Bc: \Bc^\perp \lra {}^\perp \Bc, 
 \]
 known as the {\em mutation} along $\Bc$, see \cite{BK}. It is defined as the composition
 \[
 \Bc^\perp \buildrel i_{\Bc^\perp}\over\lra \Ac \buildrel  i_{{}^\perp\Bc}^*\over\lra {}^\perp\Bc
 \]
 of the embedding of $\Bc^\perp$ and of the projection onto ${}^\perp\Bc$ along $\Bc$. 

\begin{prop}\label{prop:**}
If $\Bc$ is admissible, then the functor ${}^*i_{\Bc^\perp}$ (the projection onto $\Bc^\perp$ along $\Bc$)
has itself a left adjoint ${}^{**}i_{\Bc^\perp}: \Bc^\perp\to\Ac$ given by ${}^{**}i_{\Bc^\perp} = i_{{}^\perp \Bc}\circ M_\Bc$. 
\end{prop}

\noindent {\sl Proof:} Let $B'\in \Bc^\perp$ and $A\in\Ac$.
The definitions of $M_\Bc$ and ${}^*i_{\Bc^\perp}$ give exact triangles
\[
\begin{gathered}
i_{{}^\perp \Bc}(M_\Bc(B')) \lra B' \lra B \lra i_{{}^\perp \Bc}(M_\Bc(B'))[1], 
\\
B_1\lra A \lra {}^*i_{\Bc^\perp}(A)\lra B_1[1]
\end{gathered}
\]
with $B, B_1\in\Bc$. These triangles give canonical  identifications
\[
\begin{gathered}
\Hom_{\Ac}( i_{{}^\perp \Bc}( M_\Bc(B')), A) \,\,\buildrel (1) \over \simeq \,\, \Hom_\Ac( i_{{}^\perp \Bc}( M_\Bc(B')), {}^*i_{\Bc^\perp}(A))\,\,\buildrel (2) \over\simeq
\\
\buildrel (2) \over\simeq \,\, \Hom_\Ac(B', {}^*i_{\Bc^\perp}(A)) \,\,\buildrel (3) \over =\,\, \Hom_{\Bc^\perp}(B', {}^*i_{\Bc^\perp}(A))
\end{gathered}
\]
with the reasons being:

\begin{enumerate}
\item[(1)] since $\Hom^\bullet_{\Ac}( i_{{}^\perp \Bc}( M_\Bc(B')), B_1)=0$; 

\item[(2)] since $\Hom^\bullet_\Ac(B, {}^*i_{\Bc^\perp}(A))=0$; 

\item[(3)] since $\Bc^\perp$ is a full subcategory in $\Ac$. 
\end{enumerate}
Combined together, the identifications (1)-(3) give the claimed adjointness. \qed

\noindent {\bf C. Spherical pairs.} 
Let $\Ec_0$ be a triangulated category and $\Ec_+, \Ec_- \subset\Ec_0$ be  a pair of 
 admissible subcategories,
so that we have 
 the diagrams of embeddings
  \[
  \Ec_- \buildrel \delta_-\over \lra \Ec_0 \buildrel\delta_+\over \longleftarrow \Ec_+,
  \quad   \Ec_- ^\perp \buildrel j_-\over\lra \Ec_0 \buildrel j_+ \over\longleftarrow  \Ec_+^\perp
  \]  
with $\delta_\pm$ having a left and a right  adjoint $\delta^*_\pm, {}^*\delta_\pm$ and $j_\pm$ having a left and
a double left adjoint ${}^* j_\pm, {}^{**}j_\pm$ by Proposition \ref{prop:**}. 

\begin{Defi} The pair of  admissible subcategories $\Ec_\pm$ is called a
  {\em spherical pair}, if:
  \begin{enumerate}

  \item[(SP1)] The compositions
  \[
  {}^*j_+ \circ j_-: \Ec_-^\perp \lra\Ec_+^\perp , \quad  {}^*j_- \circ *j_+: \Ec_+^\perp \lra\Ec_-^\perp
     \]
  are equivalences. 
  
  \item[(SP2)] The compositions
  \[
  \delta^*_+ \circ \delta_-:  \Ec_-\lra \Ec_+, \quad \delta^*_- \circ \delta_+:  \Ec_+\lra  \Ec_-
     \]
  are equivalences. 
  \end{enumerate}
\end{Defi} 

\begin{rem} Note that unlike a spherical functor, the definition of a spherical pair does not
appeal to any enhancement of the triangulated category $\Ec_0$, as no functorial cones are taken. 
In the following (as elsewhere in the paper) we will, however, assume that we are
 in an enhanced situation as described in 
Appendix. 
\end{rem}

\vskip .3cm

\noindent {\bf D. From a spherical pair to a spherical functor.} Let $\Ec_\pm$ be a spherical pair.
Consider the diagram
\[
\Dc_0 = \Ec_-   \buildrel S\over\lra \,\,\,\, \Ec_+^\perp  = \Dc_1, \quad S= {}^*j_+ \circ \delta_-. 
\]

\begin{prop}
$S$ is a spherical functor. 
\end{prop}
This follows from Theorem \ref{thm:AL} (which gives that (SF1) and (SF3) imply sphericity) and from the next more
precise statement.

\begin{prop}
 \begin{enumerate}
\item[(a)] The functor $S$ has both right and left adjoints
\[
R=S^* = \delta_-^* \circ j_+, \quad L={}^*S = {}^*\delta_- \circ {}^{**}j_+. 
\]

\item[(b)] The functor 
\[
T_1 = \on{Cone}\{SR\Rightarrow\Id_{\Dc_1} \}: \,\,\Dc_1 = \Ec^\perp_+\,\, \lra \,\, \Ec_+^\perp = \Dc_1
\]
is identified with the composition ${}^* j_+ \circ j_- \circ{}^* j_- \circ j_+: \Ec^\perp_+\to\Ec^\perp_+$. 
In particular, it is invertible. 

\item[(c)] The functor
\[
T_0 = \on{Cone}\{  \Id_{\Dc_0}  \Rightarrow   RS \} [-1]: \,\,\Dc_0 = \Ec_-\,\, \lra \,\, \Ec_- = \Dc_0
\]
is identified with the composition $\delta^*_- \circ\delta_+ \circ \delta^*_+ \circ\delta_-:  \Ec_-\to\Ec_-$. 
In particular, it is invertible.

\end{enumerate}
\end{prop}

\noindent {\sl Proof:} (a)  obvious from the assumptions and the fact that the adjoint of the composition of
two functors is the composition of the adjoints in the opposite order. 
\vskip .2cm

(b) The functorial exact triangle for the semi-orthogonal decomposition $\Ec_0 = \langle \Ec_-, \Ec_-^\perp\rangle$
can be written, in our notation, as
\[
\delta_- \circ\delta_-^* \buildrel u_- \over\Longrightarrow \Id_{\Ec_0} 
\buildrel v_-\over\Longrightarrow j_-\circ {}^*j_- 
\buildrel w_-\over\Longrightarrow \delta_-\circ\delta_-^*[1]. 
\]
We note that $SR={}^*j_+\circ\delta_- \circ \delta_-^* \circ j_+$ and 
the natural transformation $c: SR\Rightarrow \Id_{\Dc_1}$ can be written as
\[
{}^*j_+ \circ_0 u_- \circ_0 j_+: \,\,\, {}^*j_+\circ\delta_-\circ \delta_-^* \circ j_+ \Longrightarrow {}^*j_+\circ j_+ \,\,=\,\,
\Id_{\Dc_1}. 
\]
Here $\circ_0$ stands for the ``$0$-composition" of a functor and a natural transformation.
This implies that
\[
\Cone(c) \,\,=\,\, {}^*j_+\circ j_- \circ {}^*j_-\circ j_+
\]
as claimed. 

\vskip .2cm

(c) We now write the functorial exact triangle for $\Ec_0=\langle \Ec_+, \Ec_+^\perp\rangle$ as
\[
\delta_+\circ\delta_+^* \buildrel u_+\over\Longrightarrow \Id_{\Ec_0} \buildrel v_+\over\Longrightarrow
j_+ \circ {}^*j_+ \buildrel w_+\over\Longrightarrow \delta_+\circ\delta_+^*[1].
\]
We note that $RS= \delta_-^*\circ j_+\circ {}^*j_+ \circ\delta_-$ and the natural transformation 
$\epsilon: \Id_{\Dc_1}\Rightarrow RS$ can be written as
\[
\delta_-^*\circ_0 v_+ \circ_0 \delta_-: \,\,\, \Id_{\Dc_0} = \delta_-^*\circ\delta_- 
\Longrightarrow \delta_-^* \circ j_+ \circ {}^*j_+\circ\delta_-.
\]
This implies that
\[
\Cone(\epsilon)[-1] \,\,=\,\, \delta_-^*\circ\delta_+ \circ\delta_+^* \circ\delta_-
\]
as claimed. \qed

\begin{ex}[(Polar coordinates)]\label{ex:conic}
Let $Y$ be a CW-complex and 
  $p: V\to Y$  a real vector bundle. Denote by $i: Y\hookrightarrow V$ the embedding
of the zero section, and by $j: V^\circ\hookrightarrow V$ the embedding of the complement of the zero section. Let $q: \Ss = V^\circ/\RR_{>0}^* \to Y$
 be the spherical bundle associated to $V$, and $\tau: V^\circ\to \Ss$ the natural projection. Let also $\pi=q\tau: V^\circ\to Y$ be the
 composite projection. 
 
 Let $\Ec_0  = \Ec_0(V)  = D^b_{\on{conic}}(V) \subset D^b(V)$ be the full subcategory of $\RR_{>0}^*$-conic complexes,
 i.e., of complexes $\Fc$ such that each $\ul H^i(\Fc)$ is constant on each orbit of  $\RR_{>0}^*$ in $V$, see
 \cite{kashiwara-schapira}.  
 Consider the subcategories
 $\Ec_\pm = \Ec_\pm(V) \subset\Ec_0$ defined as follows:
 \[
 \Ec_+\,\, =\,\,  i_* D^b(Y)\,\, \simeq \,\,D^b(Y), \quad \Ec_- \,\,=\,\, p^{-1} D^b(Y) 
 \,\,\simeq\,\,  D^b(Y). 
 \]
 Then
 \[
 \Ec_+^\perp \,\,=\,\, Rj_* \tau^{-1} D^b(\Ss) \,\,\simeq\,\,  D^b(\Ss), 
 \]
 and the functor 
 \[
 S= {}^*j_+\circ\delta_-: \Ec_-\lra\Ec_+^\perp
 \]
 is identified with the spherical functor $q^{-1}$ from Example \ref{ex:sph-p1}. This is an instance of
 the above proposition, since 
 $(\Ec_\pm)$ form a spherical
 pair. To see this, it suffices to show that both $\Ec_+$ and $\Ec_-$ are admissible
  (the conditions (SP1-2) are proved similarly to the fact that $S$ is spherical). For $\Ec_+$ this is clear
 from the standard  ``recollement'' data for complexes of sheaves on the open set $V^\circ$ and the closed set $i(Y)$.
 In particular,
 \[
 {}^\perp\Ec_+ \,\,=\,\, Rj_! \tau^{-1} D^b(\Ss). 
 \]
 For $\Ec_-$ this follows from the next remark.

  \end{ex}

\begin{rem}
In the situation of Example \ref {ex:conic}, let 
$V^*$ be the  vector bundle dual to $V$, so that   the  Fourier-Sato transform  \cite{kashiwara-schapira}
 gives an identification
 \[
  \Ec_0(V) = D^b_{\on{conic}}(V)\buildrel F\over \lra  D^b_{\on{conic}}(V^*)  = \Ec_0(V^*). 
 \]
 This identification takes the category $\Ec_\pm(V) $  to $\Ec_\mp(V^*)$.   
 \end{rem}

 \vfill\eject
 
 \section{Derived categories on the $G/P$ and Schobers on symmetric products $\hen/W$.}\label{sec:G/B}
 
 \noindent {\bf A. The braid group action on $D^b(G/B)$ and $D_\coh(T^*(G/B))$.} 
 Let $\gen$ be a split reductive Lie algebra over $\CC$, with Cartan subalgebra $\hen$ and Weyl group $W$. 
 We denote by $\hen_\RR\subset\gen_\RR$ the real parts of $\hen$ and $\gen$. 
 Let $\ds\subset \Delta_+\subset\Delta$ be the systems of  simple and positive roots of $\gen$ inside the set of all roots. 
 The complex vector space $\hen$ has an arrangement of hyperplanes $\{\alpha^\perp\}_{\alpha\in\Delta_+}$
  and so has a natural stratification by {\em flats} of this arrangements,
 see \cite{KS}.

 We consider the complex manifold $\hen/W$.  It has the induced stratification, denote it $\Sc$.
 The open stratum $(\hen/W)_0$, is the classifying space  of the {\em braid group of }  $\gen$,
 which we denote by $ \on{Br} =\on{Br} (\gen)$. By definition, $\Br$ is generated by elements $s_\alpha, \alpha\in\ds$
 subject to the relations defining the Weyl group $W$ with the exception of the relations $s_\alpha^2=1$. 
 
 \begin{ex} In the case $\gen = \gen\len_n$, the manifold $\hen/W$ is
 space of monic polynomials $f(z)$  of degree
 $n$ in one variable. The stratification $\Sc$ is the stratification by types of coincidence of roots of $f(z)$ and is labelled
 by (unordered) partitions of $n$.  The group $\on{Br}(\gen\len_n)$ is the usual Artin braid group $\Br_n$. 
  \end{ex}

 Let $B\subset G$ be the standard Borel subgroup corresponding to the choice of $\Delta_+$.
 For any $I\subset\ds$ let $P_I\supset B$ be the standard parabolic subgroup generated by $B$
 and  the exponents of the Chevalley generators associated to $(-\alpha), \alpha\notin I$. Thus $P_\emptyset =G$
 and $P_\ds=B$. 
 We write $P_\alpha=P_{\ds-\{\alpha\}}$, $\alpha\in\ds$, for the ``next to minimal" parabolic subgroup associated to $\{\alpha\}$.
 We have the $\PP^1$-fibration
 \[
 G/B \buildrel q_\alpha\over\lra G/P_\alpha
 \]
 and therefore, by Proposition \ref{prop:p1-sph}, a spherical functor
 \[
 D^b(G/P_\alpha)  \buildrel S_\alpha = q_\alpha ^{-1} \over\lra  D^b(G/B)  
 \]
 and the corresponding twist functor $T_\alpha: D^b(G/B)\to D^b(G/B)$
   written  directly as
 \be\label{eq:t-a}
 \begin{gathered}
 T_\alpha(\Fc)  = R(p_{\alpha, 2 })_!  p_{1,\alpha}^* \Fc [-1], \\
 p_{\alpha, i}: \left( \bigl( (G/B)\times_{G/P_\alpha} (G/B)\bigr) - \text{ rel. diagonal } \right)  \lra G/B, \,\, i=1,2. 
 \end{gathered}
 \ee
 
 The following fact was known  before the concept of a spherical functor was discovered, cf. \cite{BBM}.

 \begin{prop}\label{prop:br-action}
 The functors $T_\alpha$ defined by \eqref {eq:t-a}, 
  are equivalences which  satisfy the relations of $\Br(\gen)$, i.e., define an action of $\Br(\gen)$
 on the derived category $D^b(G/B)$. \qed
 \end{prop}
 
 \begin{rems} (a) Proposition \ref{prop:br-action} for $\gen =\gen\len(n)$
  has a quaternionic analog. We consider the space $\HH F_n$
 of complete flags of (left) quaternionic subspaces $V_1 \subset \cdots \subset V_{n-1}\subset\HH^n$, 
 $\dim_\HH (V_i)=i$. It fits into the  $\HH \PP^1$-fibrations
 \[
 q_i: \HH F_n \lra \HH F_n^{(i)} \,\,=\,\,\bigl\{ (V_1\subset \cdots \subset V_{i-1}\subset V_{i+1} \subset \cdots \subset \HH^n\bigr\}
 \]
 over the spaces of next-to-complete flags.  We note that $\HH \PP^1=\Ss^4$ is the 4-sphere.
 By Proposition   \ref{prop:p1-sph} this gives,
 for each $i=1, \cdots, n-1$,  a spherical functor 
 $S_i^\HH=q_i^{-1}: D^b(\HH F_n^{(i)})\to D^b(\HH F_n)$ and
 the induced twist automorphism $T_i^\HH$ of $D^b (\HH F_n)$. One can see directly by analyzing the Schubert correspondences
 in $\HH F_n$,  that the $T_i^\HH$ satisfy the relations of the Artin braid group $\Br_n$. 
 
 (b) One has also a real analog for any $\gen$, using the real loci of the flag varieties which fit into $\RR\PP^1 = \Ss^1$-fibrations. 
 
 \end{rems}
 
 A quasi-classical analog of the action in Proposition \ref{prop:br-action}  has been constructed in \cite{KhT, BR}. 
 It uses the diagram
 \[
Y_\alpha :=  T^*(G/P_\alpha)
\buildrel \rho_\alpha\over  \lla D_\alpha := T^*(G/P_\alpha) \times_{G/P_\alpha} G/B
\buildrel i_\alpha\over  \lra T^*(G/B) := X
 \]
 which produces a spherical functor
  \[
S_\alpha^\coh=  i_{\alpha *} \rho_\alpha^*: D^b_\coh T^*(G/P_\alpha) \to D^b_\coh T^*(G/B)
 \]
 It was proved in \cite{KhT} for $\gen=\gen\len_n$ and in \cite{BR} for arbitrary $\gen$,
 that the corresponding twist functors $T_\alpha^\coh$ define an action of
 $\Br(\gen)$ on $D^b_\coh T^*(G/B)$ (which, in fact, extends to an action of the affine braid group).

 These constructions  can be seen as giving {\em local systems of triangulated categories} on $(\hen/W)_0 = K(\Br(\gen), 1)$
 with  general stalk being $D^b(G/B)$, resp. $D^b_\coh(T^*(G/B))$. We would like to suggest that these local systems
 extend to natural perverse Schobers on the entire $\hen/W$. For this, we review some features of usual perverse
 sheaves in this situation. 
   
 \vskip .3cm
 
 \noindent {\bf B. Perverse sheaves on $\hen/W$ and double cubical diagrams.}
 Denote
 $\Perv(\hen/W)$ the category of perverse sheaves on $\hen/W$ smooth with respect to the stratification $\Sc$ from n$^\circ$A. 
  A complete quiver description on $\Perv(\hen/W)$ is not yet available. However, the results of
  \cite{KS} provide the following partial picture which aligns with the examples we considered earlier. 
  
  \vskip .2cm
  
  Let $\hen_\RR\subset \hen$ be the real form of $\hen$. 
  Inside $\hen/W$ we consider the ``real skeleton"  (or ``cut'') 
  \[
  K = \hen_\RR/W \,\,\,\subset \,\,\, \hen/W. 
  \]
   It can be thought of as a ``curvilinear cone",
  the image of the dominant Weyl chamber  $\hen_\RR^+\subset \hen_\RR$. Let $r=|\ds|$ be the rank of $\gen$. 
  Let $\Sc_\RR$ be the stratification of $K$ by the $2^r$ strata which are the images of the faces of $\hen_\RR^+$. 
  We denote these strata by $S_I, I\subset \ds$, so that $\dim_\RR (S_I) = |I|$. 
  
  Note that each $S_I$ is
  contractible. Therefore (cf. \cite{gelfand-macpherson})  a sheaf $\Gc$ on $K$ constructible with respect to $\Sc_\RR$ can be recovered
  from the cubical diagram of vector spaces $\Gc_I = \Gamma(S_I, \Gc)$ and {\em generalization maps}
  $\gamma_{II'}: \Gc_I\to \Gc_{I'}$, $I\subset I'$ which are required to form a representation of the poset ${\bf 2}^\ds$
  (so that the diagram is commutative). 
  
  The following is deduced by pulling $\Fc$ back to a perverse sheaf on $\hen$ smooth with
  respect to an arrangement of hyperplanes $\{\alpha^\perp\}_{\alpha\in\Delta_+}$ and applying the results of   \cite{KS}.  
  
  \begin{prop} (a) For $\Fc\in\Perv(\hen/W)$ we have $\ul\HH^i_K(\Fc)=0$ for $i\neq\dim_\RR(K)$, and the sheaf
    $\Rc_\Fc = \ul\HH^{\dim(K)}_K(\Fc)$ on $K$  is  constructible with respect to the stratification $\Sc_\RR$. 
  
  (b) Let $I\subset\ds$. Denoting $E_I(\Fc) = \Gamma(S_I, \Rc_\Fc)$, we have that $E_I: \Perv(\hen/W)\to\Vect_\k$ is an exact functor
  which takes the Verdier duality to the vector space duality. \qed
  
  \end{prop}
  
   Therefore we can associate to any  $\Fc\in\Perv(\hen/W)$  a {\em double cubical diagram} $E(\Fc)$
   formed by the vector spaces  $E_I(\Fc)$ and the maps
   \be\label{eq:EF}
   \xymatrix{
   E_I(\Fc) \ar@<.5ex>[r]^{\,\,\,\,\gamma_{IJ}} &
   E_{J}(\Fc)  \ar@<.5ex>[l]^{\,\,\,\,\,\,\delta_{JI}}
   }, \quad I\subset J. 
   \ee
  Here $\gamma_{IJ}$ is the generalization map for $\Ec_\Fc$, and $\delta_{JI}$ is the dual to the generalization map
   for $\Ec_{\Fc^*}$, where $\Fc^*$ is the perverse sheaf  Verdier dual to $\Fc$. Each of the collections $(\gamma_{IJ})$,
   $(\delta_{JI})$ forms a commutative cube. It seems very plausible that the functor $\Fc\mapsto E(\Fc)$
   from $\Perv(\hen/W)$ to the category of double cubical diagrams is fully faithful, i.e., $\Fc$ can be recovered from
   $E(\Fc)$. 
   
   \begin{ex}
   For $\gen=\sen\len_2$ we have $\hen/W=\CC$ and $K=\RR_{\geq 0}$ so the above reduces,
   very precisely,  to the construction of \S \ref{sec:disk}A, except 
   with  the disk replaced by $\CC$.

   \end{ex}
   
   \vskip .3cm
   
   \noindent {\bf C. A double cubic diagram related to flag varieties.} 
   We now note that geometry of flag varieties provides a natural double cubic diagram of categories of the
   same shape as \eqref{eq:EF}. More precisely, we have an ordinary cubical diagram of algebraic varieties $G/P_I$ and projections
   $q_{IJ}: G/P_J \to G/P_I$, $I\subset J$. We have then the double cubical diagram formed by triangulated categories
   $D^b(G/P_I)$ and the adjoint pairs of functors
   \be\label{eq:der-flag}
   \xymatrix{
   D^b(G/P_I) \ar@<.5ex>[r]^{\,\,\,\, q_{IJ}^{-1}} &
   D^b(G/P_J)  \ar@<.5ex>[l]^{\,\,\,\,\,\, Rq_{IJ *} }
   }, \quad I\subset J. 
   \ee
   As in Example \ref{ex:sph-p1}, we can also consider the left adjoint to $q_{IJ}^{-1}$ which differs from
   the right adjoint by a shift. 
   
    At the ``quasi-classical level'' one has a similar diagram formed by the categories $D^b_\coh(T^*(G/P_I))$ and the functors between them obtained
   by translating $q_{IJ}^{-1}$ and $Rq_{IJ*}$ into the language of $\Dc$-modules and passing to the associated graded modules.
   
   \vskip .2cm
   
   We would like to suggest that the diagram \eqref{eq:der-flag} comes from a more fundamental object: a perverse Schober on $\hen/W$
   extending the local system discussed in n$^\circ$A. 
   
   Similarly for the 
    quasi-classical analog with the $D^b_\coh(T^*(G/P_I))$. In this situation we have in fact more: 
   an   actions the affine braid group. 
     These  can    possibly come from  Schobers not on $\hen/W$
    but on $T/W$ where $T$ is the maximal torus in the algebraic group corresponding to $\gen$.        
      \vfill\eject

 \section{``Fukaya-style'' approach to perverse sheaves: cuts, real skeletons and Langangian varieties}\label{sec:fuk}
 
 \noindent{\bf A. Maximally real cuts.}
  Since our preliminary definitions of a perverse Schobers were based on
 quiver descriptions of perverse sheaves, let us look at some general features of such descriptions.

 Let $(X, \Sc)$ be a stratified complex manifold of dimension $n$. 
 Obtaining a quiver description of $\Perv(X, \Sc)$ requires, in particular, construction of many exact functors
$\Perv(X,\Sc)\to\Vect_\k$. Indeed, any component of the putative quiver must be such a functor. Arrows of the quiver
are then natural transformations between these exact  functors. 

The common tool  for that, used in examples in this paper, is a choice of a closed
subset $K\subset X$ with the following property: 

\begin{itemize}
\item[(Cut)] For any $\Fc\in\Perv(X,\Sc)$ we have $\ul\HH^i_K(\Fc) = 0$, $i\neq n$.  
\end{itemize}
This property implies that the functor of abelian categories
\[
\Rc: \Fc\longmapsto \Rc(\Fc) :=  \ul\HH^n_K(\Fc), \quad \Perv(X,\Sc)\lra \on{Sh}_K
\]
is exact. So stalk of $\Rc(\Fc)$ at any point $x\in K$ can be used as a component of a quiver,
while generalization maps between the stalks provide some of the arrows. 

We will call each $K$ satisfying (Cut) an  {\em (admissible) cut} for $(X, \Sc)$ and denote by $\Cc(X,\Sc)$
the set of all such cuts. It is natural to look for a description of $\Perv(X,\Sc)$ in terms of some data
associated to all the cuts.

 Recall that by the  Riemann-Hilbert correspondence we can realize each  $\Fc\in\Perv(X,\Sc)$ as $\ul{R\Hom}_{\Dc_X}(\Mc, \Oc_X)$
for a  left $\Dc_X$-module $\Mc$ whose characteristic variety satisfies the inclusion
 \[
\Ch(\Mc)  \,\, \subset  \,\, \Lambda_\Sc :=  \bigcup _{\alpha} \ol{ T^*_{X_\alpha} X}. 
\]
 We denote by $_\Dc\Mod(X,\Sc)$ the category of such $\Dc_X$-modules. 
In these terms, a more detailed scenario (sufficient condition)  for (Cut) to hold would be  for $K$ to satisfy the following two properties:

\begin{itemize}
\item[(Cut1)] We have $\ul\HH^i_K(\Oc_X)=0$ for $i\neq n$. 

\item[(Cut2)] Assuming (Cut1), 
the sheaf $\Bc_K= \ul\HH^n_K(\Oc_X)$ (which, considered as a sheaf on $X$,  is automatically a sheaf of left $\Dc_X$-modules),
satisfies 
\[
\ul{\on{Ext}}^j_{\Dc_X}(\Mc, \Bc_K) = 0, \quad \forall \Mc \,\, \in\,\, {}_\Dc\Mod(X,\Sc), \,\,\,  j>0. 
\]
\end{itemize}

The property (Cut1) is not related to  a choice of $\Sc$ and holds for any {\em totally real subset} of $X$, see \cite{harvey}.
More precisely, we recall, see  \cite{baouendi} for background:

\begin{Defi}
(a) Let $V$ be a  $\CC$-vector space of finite  dimension $n$. A real  subspace $L\subset V$ is called
{\em totally real}, if $L\cap iL = 0$. We say that $L$ is {\em maximally real}, if it is totally real of dimension $n$.  

(b) Let $X$ be a complex manifold of dimension $n$. A $C^\infty$-submanifold $K\subset X$   is called
{\em totally real}, resp. {\em maximally real} if for each $x\in K$ the subspace $T_xK\subset T_xX$ is totally real resp. maximally real.

\end{Defi}

The results of Harvey \cite{harvey} imply:

\begin{prop}
 Any closed subset of a totally real submanifold of $X$ satisfies (Cut1). 
\end{prop}

More precisely, the result of {\em loc. cit.} is for arbitrary  totally real {\em subsets} of $X$, a class of sets which includes
totally real submanifolds \cite[\S 3.6 Ex. 2]{harvey} and is closed under passing to closed subsets \cite[Cor. 3.2]{harvey}. 

For example, $\RR^n$, as well as $\RR_{\geq 0}^n$ satisfies (Cut1). The sheaf $\Bc_{\RR^n}$ is the 
sheaf of hyperfunctions of Sato \cite{SKK}. 

\vskip .2cm

 We now consider the condition (Cut2). Sufficient criteria for it to hold
 were given by Lebeau \cite{lebeau} and Honda-Schapira \cite{honda-schapira}.
 The criterion of \cite{honda-schapira} is based on the concept of {\em positive position} of  two real Lagrangian submanifolds
 in a complex symplectic manifold such as $T^*X$. We do not recall this concept
 here,  referring to    \cite{schapira-pos} and references therein for more background. Informally, the
 essense of the criterion can be formulated like this.
 \be
 \begin{gathered}
 \text{ \em For (Cut2) to hold, $K$ must be ``maximally real with respect to} \\
 \text{\em the stratification $\Sc$'', in particular,
 the intersection of $K$} \\
 \text{\em with each stratum $X_\alpha$ should be maximally real in $X_\alpha$. 
 }
 \end{gathered}
 \ee
 Here is a precise but more restrictive statement which is a reformulation of 
   Example 1 of  \cite{honda-schapira}.
 
 \begin{prop} Suppose $\Sc$ is such that $\Lambda_\Sc$ is contained in the
   union of $T^*_{Z_\beta} X$ for a collection of {\em smooth closed}
 complex  submanifolds $Z_\beta \subset X$. Suppose $K\subset X$ is a maximally real analytic
 submanifold such that  each $K\cap Z_\beta$ is maximally real in $Z_\beta$. Then $K$ satisfies
 (Cut2). \qed
 \end{prop}
 
 \vskip .3cm
 
 \noindent {\bf B. Maximally real vs. Lagrangian cuts.} 
 Let $V$ be a complex vector space of dimension $n$ and $G_\RR(n, V)$ the Grassmannian of real $n$-dimensional subspaces
 in $V$. We denote by $G^{\max} (V)$ the open subset in $G_\RR(n, V)$ formed by maximally real subspaces. 
 
 Suppose $V$ is equipped with a positive definite hermitian form $h$. Separating the real and imaginary parts $h= g+i \omega$, 
 we have that $\omega$ is a symplectic form on $V$. Let $\LG_\omega (V) $ be the closed subset in
 $G_\RR(n,V)$ formed by
 subspaces Lagrangian with respect to $\omega$. The following is well known.
 
 \begin{prop}
 $\LG_\omega(V)$ is contained in $G^\max(V)$ and the embedding is a homotopy equivalence.
  In other words, $\LG_\omega(V)$ can be seen as a compact form of $G^\max(V)$. 
 \end{prop}
 
\noindent {\sl Proof:}  For $V=\CC^n$ with the standard hermitian form we have
 \[
 \LG_\omega(V) \,\,=\,\, O(2n)/U(n) \,\,\subset \,\, GL(2n, \RR)/GL(n, \CC) \,\,=\,\, G^\max(V). 
 \]
 \qed

 Let now $X$ be a complex manifold equipped with a K\"ahler metric $h = g+i\omega$. Then $\omega$
 makes $X$ into a symplectic manifold, and we have
 
 \begin{cor}
 Any Lagrangian submanifold of $X$ is maximally real. \qed
 \end{cor}

  Note that all the cuts used in this paper as well as in \cite{GGM} \cite{KS}, are Lagrangian. 
  This suggests a possibility of describing more general $\Perv(X,\Sc)$ in terms of data coming from
  Lagrangian cuts.
  
  \vskip .2cm
  
  On the other hand, Lagrangian submanifolds of a K\"ahler manifold $X$ are organized into a far-reaching
  structure: the {\em Fukaya category} $\FF(X)$ (as well as its modifications such as the Fukaya-Seidel and wrapped Fukaya categories),
  see \cite{FOOO} \cite{FSbook} for background. We briefly recall that $\FF(X)$ is an ({\em a priori} partially defined)
  $\CC$-linear
  $A_\infty$ category whose objects are, in the simplest setting, compact Lagrangian submanifolds $K\subset X$. 
  The space $\Hom(K_1, K_2)$ is defined when $K_1$ and $K_2$ meet transversely and in this case is formally spanned
  by the set $K_1\cap K_2$ (with appropriate grading, see {\em loc. cit.}). 
  The $A_\infty$-composition
  \[
  \mu_{K_1, \cdots, K_m}: \,\,\, \bigotimes_{i=1}^{m-1} \Hom(K_i, K_{i+1}) \lra\Hom(K_1, K_m), \quad m\geq 1, 
  \]
  is given by counting holomorphic\footnote{
  Here we leave aside the additional complication that it may be necessary to deform the complex structure on $X$
  to ensure generic position. 
  }
  disks $D\subset X$ with boundary on $K_1\cup \cdots\cup K_m$. 
  Each $D$ gives a contribution $\exp\bigl(-\int_D \omega\bigr)$ to the appropriate matrix element of $\mu_{K_1, \cdots, K_m}$.

 \vskip .2cm
 
 One can expect that the Fukaya structure on the collection of admissible Lagrangian cuts for $(X,\Sc)$ has some
 significance for explicit description of $\Perv(X,\Sc)$. In particular, a holomorphic disk $D$ with boundary on the union of
 admissible cuts $K_1, \cdots, K_m$ may provide a link between the sheaves $\ul\HH^n_{K_\nu}(\Fc)$,
 $\Fc\in \Perv(X,\Sc)$ for different $\nu$, via the structure of perverse sheaves on a disk (\S \ref{sec:disk-several}). 
 
 \begin{rems}
 (a) Any collection $K_1, \cdots, K_m$ of maximally real (not necessarily Lagrangian)
 submanifolds in $X$ provides a natural boundary condition for holomorphic disks. It has been noticed
 \cite{FSbook} that many ingredients of the Fukaya category construction have an ``intrinsic" meaning and
 can be defined without explciit reference to the symplectic structure. In particular:
 
 \begin{itemize}
 \item The role of the grading on $\Hom(K_1, K_2)$ is to make the count of disks possible by ensuring that
 in evaluating matrix elements of $\mu$ {\em between basis vectors of the same degree},
 the index of the linearlized elliptic problem is $0$ (so the corresponding linear $\ol\partial$-operator
 is, generically, invertible).
 
 \item The quantity  $\exp\bigl(-\int_D \omega\bigr)$ can be intepreted as the determinant of the invertible
 $\ol\partial$-operator above. 
 
 \end{itemize}
 
 \noindent They could therefore  make sense  for more general maximally real submanifolds. The usual technical reason for 
 restricting to
 Lagrangian submanifolds in defining $\FF(X)$ is that the Gromov compactness theorem
 (which, via the properties of 1-dimensional moduli spaces, is used to prove the $A_\infty$-axioms)
   has originally been
 established only in that setting. However, see \cite{frauenfelder} for a recent generalization to the maximally real case. 
 
 \vskip .2cm
 
 (b) the above mentioned role of cuts (in particular, Lagrangian cuts) for description of perverse sheaves on $X$,
 is different from the classical ``microlocal'' point of view \cite {kashiwara-schapira} which empasizes complex,
 conic Lagrangian subvarieties in $T^*X$ rather than real ones in $X$ itself. For the relation of the microlocal
 approach to the Fukaya category (of $T^*X$), see \cite{nadler-zaslow}\cite{nadler}. 
 A more general idea that the Fukaya category of any symplectic manifold should have an interpretation in terms of its
   geometric quantization, was proposed  earlier in \cite{BS}.

 \end{rems}

 \vskip .3cm
 
 \noindent {\bf C. Perverse Schobers as coefficients data for Fukaya categories.} The Fukaya category 
 $\FF(X)$ of a symplectic manifold $X$ can be seen as a categorification of its middle (co)homology (or, rather,
 of the part represented by Lagrangian cycles). Here ``cohomology'' is understood as $H^\bullet(X, \CC)$, the cohomology with
 constant coefficients, a particular case of a more general concept of $H^\bullet(X, \Fc)$, the cohomology with coefficients in a sheaf $\Fc$
 or  even more generally, in a complex of sheaves. 
 
 \vskip .2cm
 
 This point of view leads to the idea of introducing coefficients into the definition of the Fukaya category as well. 
  It  was proposed by M. Kontsevich with the goal of understanding the  (usual) Fukaya category of a manifold by projecting it
 onto a manifold  of smaller dimension.
 
 \vskip .2cm
 
  In this direction we would like to suggest that perverse Schobers are the right coefficient data for defining Fukaya categories.
  That is, to a perverse Schober $\Sen$ on a K\"ahler manifold $X$ there should be naturally associated a triangulated
  category $\FF(X,\Sen)$, which for the constant Schober  ${D^b(\Vect_\CC)}$ reduces to $\FF(X)$.
    If we think about perverse sheaves in terms of some vector space data associated to Lagrangian cuts and
  then categorify these data to define Schobers, then it is natural to try to define $\FF(X,\Sen)$ in terms of these
  categorified data.
  
  \vskip .2cm
  
  For example,   when $X$ is a Riemann surface with marked points, the Fukaya category of $X$
  with coefficients in a constant ($\ZZ/2$-graded) triangulated category $\Ac$ was defined in \cite{DK}.
  One can easily modify this definition  when $\Ac$ is replaced by a local system of
  triangulated categories on $X$, i.e.,  by the next simplest instance of a perverse Schober. 
  We leave the case of an arbitrary perverse Schober on a Riemann surface 
 (\S \ref{sec:disk-several}E) 
  for
  future work.

 \vfill\eject

\appendix

\section {Conventions.} 
 For the considerations of this paper to make sense,  we should work in some framework of
  ``refined'', or ``enhanced''  triangulated categores having functorial cones. Let us describe one such framework, to be followed in the 
  main body. 
  (Another  one  would be that of  stable $\infty$-categories
\cite{lurie-stable}.)

\vskip .2cm

 All our categories will be $\k$-linear, where $\k$ is a fixed base field. 
We recall that Tabuada \cite{Ta1, Ta2} has introduced a {\em Morita model structure} on the category of
dg-categories. Dg-categories fibrant with respect to this model structure are called {\em perfect}, see also \cite{TV-moduli}. 
Perfect dg-categories form a subclass of pre-triangulated categories in the sense of  \cite{BK-enhanced}.
That is,   a perfect dg-category $\Ac$ gives rise to a triangulated category $H^0(\Ac)$ which is, in addition, 
idempotent complete.  

\vskip .2cm

In the main body of this paper the word ``triangulated category" will always mean ``a triangulated category $\Dc$ together with
an identification $\Dc\simeq H^0(\Ac)$ where $\Ac$ is a perfect dg-category", and ``exact functor" will mean
``an exact functor of triangulated categories obtained, by passing to $H^0$, from a dg-functor between perfect dg-categories"
and similarly for ``natural transformation".
With this understaning we will speak about  ``the'' exact functor 
$\on{Cone} \{T: F\Rightarrow G\}$ where $T$ is a natural transformation
of exact functors. 

\vskip .2cm

For a finite CW-complex $Z$ we denote $D^b(Z)$ the bounded derived category of all sheaves of $\k$-vector spaces on $Z$.
For a smooth complex algebraic variety $X$ we denote by $D^b_\coh(X)$ the bounded derived category of
coherent sheaves on $X$. 

 The condition of perversity for constructible complexes on a complex manifold $X$ is normalized in such a way that
 a constant sheaf in degree $0$ is perverse.

\vfill\eject

  \let\thefootnote\relax\footnote {
M.K.: Kavli Institute for Physics and Mathematics of the Universe (WPI), 5-1-5 Kashiwanoha, Kashiwa-shi, Chiba, 277-8583, Japan.
Email: {mikhail.kapranov@ipmu.jp}

\vskip .2cm

V.S.: Institut de Math\'ematiques de Toulouse, Universit\'e Paul Sabatier, 118 route de Narbonne, 
31062 Toulouse, France, schechtman@math.ups-toulouse.fr 

}


\begin{thebibliography}{100}

\small

 
  
 \bibitem[A1]{A1}  R. Anno Spherical functors, arXiv:0711.4409. 
 
 \bibitem[A2]{A2} R. Anno. 
Affine tangles and irreducible exotic sheaves, arXiv:0802.1070. 
  
\bibitem[AL1]{AL1} R. Anno, T. Logvinenko.  Orthogonally spherical objects 
and spherical fibrations, arXiv:1011.0707. 

\bibitem[AL2]{AL2}     R. Anno, T. Logvinenko.  Spherical DG-functors, 
arXiv:1309.5035.

\bibitem[BER]{baouendi} M.~Salah Bauendi, P. Ebenfeld, L.~P.~Rotschild.
\newblock Real Submanifolds in Complex Space and Their Embeddings.
\newblock Princeton Univ. Press, 1999. 

\bibitem[Be]{beil-gluing} A.~Beilinson.
\newblock How to glue perverse sheaves. In: $K$-theory, arithmetic and geometry (Moscow, 1984), 
\newblock {\em  Lecture Notes in Math.} {\bf 1289}, 
Springer-Verlag, 1987, 42 - 51.

\bibitem[BBD]{BBD}  A.~Beilinson, I.~Bernstein, P.~Deligne.
\newblock Faisceaux Pervers, {\em Ast\'erisque} {\bf 100}, 1982. 

\bibitem[BBM]{BBM} A.~Beilinson,  R.~Bezrukavnikov, I.~Mirkovic.
\newblock Tilting exercises. \newblock {\em  Mosc. Math. J.}
 {\bf 4} (2004),  547–557. 
 
 
 \bibitem[BR]{BR}
 R.~Bezrukavnikov, S.~ Riche.
  \newblock Affine braid group actions on derived categories of Springer resolutions.
  \newblock {\em  Ann. Sci. \'Ec. Norm. Sup\'er.} {\bf 45}  (2012), 535-599. 
  
  
  \bibitem[BK1]{BK} A. I. Bondal, M. M. Kapranov.
  Representable functors, Serre functors and mutations.
   {\em Math. USSR Izv. } {\bf 35}  (1990), 519-541. 
   
   \bibitem[BK2]{BK-enhanced} A. I. Bondal, M. M. Kapranov. Enhanced triangulated categories.
   {\em Math. USSR Sbornik},  {\bf 70} (1991) 93-107. 
   
   \bibitem[BS]{BS} P.~Bressler, Y.~Soibelman.
  Mirror symmetry and deformation quantization. 
  \newblock arXiv hep-th/0202128. 
  
     
 \bibitem[BT]{BT}  C.Brav, H.Thomas, Braid groups and Kleinian singularities, 
{\em Math. Ann.} {\bf 351} (2011), 1005 - 1017. 

\bibitem[CG]{CG} Neil Chriss, Victor Ginzburg, Representation Theory and 
Complex Geometry.  Birkh\"auser, 1997.

\bibitem[Cu]{curry} J.~Curry.
\newblock Sheaves, Cosheaves and Applications.
\newblock 
arXiv:1303.3255.

\bibitem[DK]{DK} T.~Dyckerhoff, M.~Kapranov.
Triangulated surfaces in triangulated categories.
ArXiv: 1306.2545. 
 

\bibitem[FZ]{frauenfelder} U.~Frauenfelder, K.~Zehmisch.
\newblock  Gromov compactness for holomorphic discs with totally real boundary conditions. 
\newblock arXiv:1403.6139.  

\bibitem[FO$^3$]{FOOO} 
K.~Fukaya, Y.-G.~Oh, H.~Ohta, K.~Ono.
\newblock Lagrangian Intersection Floer Theory.
\newblock Amer. Math. Soc. Publ. 2009. 

\bibitem[GGM]{GGM} A.~Galligo, M.~Granger, P.~ Maisonobe.
\newblock $\Dc$-modules et faisceaux pervers dont le support singulier
est un croisement normal. \newblock {\em Ann. Inst. Fourier Grenoble},
{\bf 35} (1985), 1-48. 
 
\bibitem[GM]{gelfand-macpherson}
S.~Gelfand, R.~D.~MacPherson.  
\newblock Verma modules and Schubert cells: a dictionary. 
\newblock in: {\em ``Paul Dubreil and Marie-Paule Malliavin Algebra Seminar", 34th Year (Paris, 1981)}, pp. 1-50.
\newblock
 Lecture Notes in Math. {\bf 924}, Springer, Berlin-New York, 1982.
 
 \bibitem[GMV]{gelfand-MV} 
 S.~Gelfand, R.~D.~MacPherson, K.~Vilonen.
 \newblock Perverse sheaves and quivers.
 \newblock {\em Duke Math. J. } {\bf 3} (1996), 621-643. 
 
 
 \bibitem[Har]{harvey} R.~Harvey. The theory of hyperfunctions on totally real subsets of a complex manifold with
 application to extension problems. {\em Amer. J. Math.} {\bf 91} (1969),  853-873. 
 
 \bibitem[HS]{honda-schapira} N.~Honda, P. ~Schapira.
\newblock
A vanishing theorem for holonomic modules with positive characteristic varieties.
{\em Publ. RIMS Kyoto Univ.} {\bf 26} (1990), 529-534. 
 
 \bibitem[KS1]{KS} M.Kapranov, V.Schechtman, Perverse sheaves over real hyperplane 
arrangements, arXiv:1403.5800, {\it Ann. of Math.} (2016), to appear. 

\bibitem[KS2]{kashiwara-schapira} M. Kashiwara, P. Schapira. Sheaves on manifolds. 
\newblock
Springer, 1990. 

\bibitem[KhT]{KhT} M. Khovanov, R. Thomas.
 Braid cobordisms, triangulated categories, and flag varieties, arXiv:math/0609335.  
 
 \bibitem[KuL]{kuznetsov} A.~Kuznetsov, V.~Lunts.
 Categorical resolution of irrational singularities, \newblock 
 arXiv:1212.6170.  
 
 
\bibitem[Le]{lebeau} G.~Lebeau. Annulation de la cohomologie hyperfonction de certains
modules holonomes. {\em C.R. Acad. Sci. Paris, S\'er. A,} t.{\bf 290} (1980), 313-316.  
 
 \bibitem[Lu]{lurie-stable} J. Lurie. Stable $\infty$-categories, arXiv math/0608228. 
 
  
 \bibitem[Na]{nadler} D.~Nadler. Fukaya categories as categorical Morse homology.
 \newblock  {\em SIGMA Symmetry Integrability Geom. Methods Appl.} {\bf 10} (2014), Paper 018. 

\bibitem[NZ]{nadler-zaslow}  
  D.~ Nadler, E.~Zaslow.
  Constructible sheaves and the Fukaya category. 
  \newblock {\em J. Amer. Math. Soc.}  {\bf 22}  (2009) 233-286. 

 
 
 
 \bibitem[SKK]{SKK} M. Sato, T. Kawai, M. Kashiwara. Hyperfunctions and pseudo-differential equations.
 {\em Lecture Notes in Math.}  {\bf 287} (1973), 265 - 529. 
 
 \bibitem[Sch]{schapira-pos} P. Schapira.
  \newblock
 Conditions de positivit\'e dans une vari\'et\'e symplectique complexe. Application a l' \'etude des microfonctions.
 \newblock{\em Ann. Sci. ENS} {\bf 14} (1991), 121-139. 
 
 
 \bibitem[Se]{FSbook} P. Seidel, {Fukaya categories and Picard-Lefschetz theory}, European Mathematical Society, Zurich, 2008.
 
\bibitem[ST]{ST} P. Seidel, R. Thomas, Braid group actions on derived categories 
of coherent sheaves, arXiv:math/0001043. 


\bibitem[Ta1]{Ta1} G.  Tabuada.
 Une structure de cat\'egorie de mod\`eles de Quillen sur la cat\'egorie des dg-cat\'egories. 
 {\em C. R. Math. Acad. Sci. Paris}  {\bf 340} (2005) 15-19. 

\bibitem[Ta2]{Ta2} G. Tabuada. Th\'eorie homotopique des DG-cat\'egories,
arXiv:0710.4303. 

\bibitem[TV]{TV-moduli} B. To\"en, M. Vaqui\'e. \newblock 
 Moduli of objects in dg-categories. {\em Ann. Sci. \'Ecole Norm. Sup. (4)} {\bf 40} (2007)  387-444. 
  

  \end{thebibliography}
\end{document}